\documentclass[11pt]{amsart}    

\usepackage{graphicx}

\usepackage{amsmath,amssymb,amsfonts}
\usepackage{color}                    
\usepackage[colorlinks,citecolor=blue]{hyperref}

\usepackage{graphics}
\usepackage{pgf,tikz}
\usetikzlibrary{arrows}
\usetikzlibrary{calc}
\definecolor{uququq}{rgb}{0.25,0.25,0.25}

\usepackage{amsmath,amssymb,amsfonts}
\usepackage{color}                    
\usepackage[colorlinks,citecolor=blue]{hyperref}
\newtheorem{theo}{Theorem} 
\newtheorem{pro} {Proposition}

\newtheorem{remark}{Remark}
\newtheorem{lem}{Lemma}

\newtheorem{definition}{Definition}
\newtheorem{cor}{Corollary}

\newcommand{\rr}{{\mathbb{R}}} 
\newcommand{\RR}{{\mathbb{R}}} \newcommand{\R}{{\mathbb{R}}}
\newcommand{\nn}{{\mathbb{N}}} \newcommand{\N}{{\mathbb{N}}}
\newcommand{\Q}{{\mathbb{Q}}}

\newcommand\dimh{\dim_{\cH}}

\newcommand{\al}{\alpha }
\newcommand{\be}{\beta }
\newcommand{\ga}{\gamma }
\newcommand\de{\delta}

\newcommand\ep{\varepsilon}\newcommand{\e}{\varepsilon}

\newcommand{\pp}{{\mathbb{P}}}
\newcommand{\E}{{\mathbb{E}}}
\newcommand{\indiq}{{\bf 1}}

\newcommand{\ds}{\displaystyle}
\newcommand{\wt}{\widetilde}
\newcommand{\wh}{\widehat}
\newcommand{\ov}{\overline}
\newcommand{\un}{\underline}

\newenvironment{preuve}{\noindent {\it Proof}}{\hfill$\square$}

\newcommand\zu{[0,1]}

\newcommand{\supp}{\text{Supp\,}}

\newcommand{\cB}{{\mathcal B}}
\newcommand{\cC}{{\mathcal C}}
\newcommand{\cD}{{\mathcal D}}
\newcommand{\cE}{{\mathcal E}}
\newcommand{\cF}{{\mathcal F}}

\newcommand{\cH}{{\mathcal H}}

\newcommand{\cJ}{{\mathcal J}}

\newcommand{\cL}{{\mathcal L}}
\newcommand{\cM}{{\mathcal M}}

\newcommand{\cO}{{\mathcal O}}
\newcommand{\cP}{{\mathcal P}}

\newcommand{\cR}{{\mathcal R}}

\newcommand{\cT}{{\mathcal T}}

\def \( {\left(}
\def \) {\right)}
\def \[ {\left\lbrack}
\def \] {\right\rbrack}
\def \lba {\left|}
\def \rba {\right|}
\def \llb {\left\lbrace}
\def \rrb {\right\rbrace}
\def \< {\left\langle}
\def \> {\right\rangle}
\def\sk{\smallskip}
\def\mk{\medskip}

\begin{document}

\title[{Stable-like occupation measure}  ]{Multifractal analysis for the  occupation measure of stable-like processes}\thanks{The first author was   partially supported by the ANR project MUTADIS, ANR-11-JS01-0009.  The second author was supported by grants from R\'egion Ile-de-France.
}


\author{St\'ephane SEURET}
\address{ Universit\'e Paris-Est ,              LAMA (UMR 8050), UPEMLV, UPEC, CNRS, F-94010, Cr\'eteil, France}
\email{seuret@u-pec.fr}   
\author{   Xiaochuan YANG}
\address     {        Universit\'e Paris-Est,
              LAMA (UMR 8050), UPEMLV, UPEC, CNRS, F-94010, Cr\'eteil, France}
\email{xiaochuan.yang@u-pec.fr}  

 
\begin{abstract}
In this article, we investigate   the local behaviors of the occupation measure $\mu$ of a class of real-valued Markov processes $\mathcal{M}$, defined via a SDE.   This (random) measure  describes the time spent in each set $A\subset \R$ by the sample paths of $\mathcal{M}$. We compute the multifractal spectrum of $\mu$, which turns out to be random, depending on the trajectory. This  remarkable property is in sharp contrast with the results  previously obtained on occupation measures of other processes (such as L\'evy processes), since the multifractal spectrum is usually  deterministic, almost surely. In addition, the shape of this multifractal spectrum is very original, reflecting  the richness and variety of the local behaviors. The proof is based on new methods, which lead for instance to fine estimates on Hausdorff dimensions of certain jump configurations in Poisson point processes.
\end{abstract}

\maketitle

\section{Introduction}
The occupation measure of a $\R^d$-valued stochastic process $(X_t)_{t\geq 0}$ describes the time spent by $X$ in any borelian set $A\subset \rr^d$. It is the natural measure supported on the range of the process $X$,  and plays an important role in describing the different fractal dimensions of the range of $X$.  Local regularity results for the occupation measure and its density when it  exists (often called local times if $X$ is Markovian) yield considerable information about the path regularity of the process itself, see the   survey article by Geman and Horowitz \cite{geman1980occupation} on this subject. 

We describe the local behaviors of this occupation measure via its multifractal analysis.    Multifractal analysis is now identified as a fruitful approach to provide organized information on the fluctuation of the local regularity of functions and measures, see for instance \cite{jaffard2004survey,falconer2003}. Its use in the study of pointwise regularity of stochastic processes and random measures has attracted much attention in recent years, e.g. (time changed) L\'evy processes \cite{jaffard1999levy,barral2007timechange,durand2009singularity,durand2012levyfields,balanca2014levy}, stochastic differential equations with jumps \cite{barral2010increasing,yang2015jumpdiffusion,xu2015boltzmann}, the branching measure on the boundary of a Galton-Watson tree \cite{morters2002thinbranching,morters2004branchingmeasure}, local times of a continuous random tree \cite{berestycki2008coalescent,balanca2015uniform}, SPDE \cite{perkins1998superbrownian,mytnik2015superstable,khoshnevisan2015multifractal}, Brownian and stable occupation measure \cite{hu1997occupation,shieh1998occupation,marsalle1999localtime,hu2000occupation,dembo2000thickspatial}, amongst many other references.

In this article, we obtain the almost-sure multifractal  spectrum of the occupation measure of stable-like jump diffusions, which turns out to be random, depending on the trajectory. This  remarkable property is in sharp contrast with the results  previously obtained on occupation measures of other processes (such as L\'evy processes), since the multifractal spectrum is usually  deterministic, almost surely. In addition, the shape of this multifractal spectrum is very original, reflecting  the richness and variety of the local behaviors. The proof is based on new methods, which lead for instance to fine estimates on Hausdorff dimensions of certain jump configurations in Poisson point processes.

\mk

We first introduce the class of processes we focus on.
 
\begin{definition}\label{defi1}
\label{beta}
Let $\ep_0>0$, and $\be: \rr\to [\ep_0,1-\ep_0]$ be a nowhere constant non-decreasing, Lipschitz continuous map.
The stable-like processes $\cM$ are pure jump Markov processes whose generator can be written as 
\begin{align}
\label{defstablelike}
\mathfrak{L}f(x)= \int_0^1 (f(x+u)-f(x)) \be(x) u^{-1-\be(x)} du.
\end{align} 
\end{definition}

 Introduced by Bass \cite{bass1988uniqueness} in the late 80's by solving a martingale problem, this class of processes has sample paths whose characteristics  change as time passes, which is a relevant feature when modeling real data   (e.g. financial, geographical data).  Roughly speaking, the stable-like processes behave locally like a stable process, but the stability parameter evolves following the current position of the process, see \cite{barral2010increasing} or \cite{yang2015jumpdiffusion} for an explanation from the tangent processes point of view. 

\sk
Let $\cM=\{\cM_t, t\in[0,1]\}$ be a stable-like process associated with a given function $x\mapsto\be(x)$ as in Definition \ref{defi1}. Our purpose is to describe the local behaviors of   the occupation measure of $\cM$ defined as 
\begin{align}
\label{defoccup}
\mu(A) = \int_0^1 \indiq_A(\cM_t)dt.
\end{align}
It depicts how long    $\cM$ stays in any Borel set $A\subset\rr$. We investigate the possible local dimensions for $\mu$, as well as its multifractal spectrum. Let us recall these notions.

\begin{definition} Let $\nu$ be a positive measure on $\rr$, whose support  is $\supp (\nu):=\{x\in \R: \forall \, r>0, \ \nu(B(x,r))>0\}$.
The upper  local dimension of $\nu$ at the point $x\in\supp(\nu)$ is defined  by
\begin{align*}
\overline{\dim}(\nu,x) = \limsup_{r\downarrow 0} \frac{\log \nu(B(x,r))}{\log r}. 
\end{align*}
Similarly, the lower local dimension  of $\nu$ at $x $ is  
\begin{align*}
\underline{\dim}(\nu,x) = \liminf_{r\downarrow 0} \frac{\log \nu(B(x,r))}{\log r}. 
\end{align*}
When $\overline{\dim}(\nu,x)$ and $\underline{\dim}(\nu,x)$  coincide at   $x$, the common value is denoted by  $\dim(\nu,x)$, the local dimension of $\nu$ at $x$.
\end{definition}

Our aim is to investigate two   multifractal spectra of the occupation measure $\mu$ associated with stable-like processes, related to these local dimensions. Let $\dim_{\cH}$ stand for the Hausdorff dimension in $\R$, with the convention that $\dim_{\cH}(\emptyset)=-\infty$. The first multifractal spectrum (in {\em space}) is defined as follows.
\begin{definition}\label{defi3} 
Let $\cO$ be an open set in $\R$, and $\nu$ a Borel measure on $\R$. Consider the level sets
\begin{align*}
  \ov E_\nu(\cO, h)&=\{ x \in \cO \cap \supp(\nu): \ov \dim(\nu,x) =h \}\\
    \un E_\nu(\cO,h)&=\{ x \in \cO \cap \supp(\nu): \un \dim(\nu,x) =h \}.
\end{align*}
 
The upper and lower multifractal spectrum of    $\nu$  are the mappings
\begin{align*}
 \overline{d}_\nu (\cO,\cdot) :  &  \ h\mapsto \dim_{\cH} \ov E_\nu(\cO, h),\\  
 \underline{d}_\nu (\cO,\cdot) : & \ h\mapsto \dim_{\cH} \un E_\nu( \cO, h).
  \end{align*}
  \end{definition}

 The famous paper by Hu and Taylor \cite{hu1997occupation} states that for every  $\alpha$-stable subordinator $\cL^\alpha$ whose occupation measure is denoted by $\mu^\alpha$, almost surely for all $x\in \supp \mu^\alpha$, 
 \begin{equation}\label{spec-stable}
 \underline{\dim}(\mu^\alpha,x) = \alpha \  \mbox{ and }   \ \overline{\dim}(\mu^\alpha,x)\in [\alpha,2\alpha].
 \end{equation}

It is a classical result \cite{blumenthal1960stable} that when $\alpha\in (0,1)$, the image of any interval $I$ by $\cL^\alpha$ has Hausdorff dimension $\alpha$, almost surely. This implies that the support of $\mu^\alpha$ has Hausdorff dimension $\alpha$, almost surely. 

With all this in mind, the lower spectrum of $\mu^\alpha$ is trivial: for each open interval $\cO$ intersecting $\supp(\mu)$,  one has $\un d_{\mu_{\al}} (\cO, h)= \begin{cases} \ \  \alpha & \mbox{ when } h =\alpha,\\  -\infty  & \mbox{ when }h\ne \alpha.\end{cases}$

Hu and Taylor also prove that the upper spectrum is much more interesting (see Figure \ref{fig-stable}):  Almost surely,  for each open interval $\cO$ that intersects  $\supp(\mu^\al)$, 
 \begin{equation}\label{spec-stable2}
 \ov d_{\mu_{\alpha}} (\cO, h)=    g_\alpha(h) :=\begin{cases}  \alpha\left(\frac{2\alpha}{h}-1\right) &\ \mbox{ when }h\in [\alpha,2\alpha ] , \\ -\infty & \mbox{ otherwise.} \end{cases}
 \end{equation} 
 
\begin{center}
\begin{figure}
\begin{center}
\begin{tikzpicture}[xscale=2.5, yscale=2.4]
\draw [fill,->] (-0.2,-0)--(1.6,0) node [right] {$h$  };
\draw [fill,->] (0,-0.2)--(0,0.89) node [above] {$\ov d_{\mu^\alpha}(h)=g_\alpha(h)$ }; 
\draw  [fill] (0,0.6)circle [radius=0.03]  node [left] {$\alpha$ \ } (0,0.6) [dotted] (0,0.6)--(0.6,0.6)  [fill] (0.6,0.6)circle [radius=0.03] ;
\draw  [fill] (0.6,0)circle [radius=0.03]  node [below] {$\alpha$ \ } (0,0.6) [dotted] (0.6,0)--(0.6,0.6);
\draw  [fill] (1.2,0) circle [radius=0.03]  node [below] {$2\alpha$ \ };
 
 \draw [thick, domain=0.6:1.2, color=black]  plot ({\x},  {0.6*(2*0.6/\x-1)});
\draw [fill] (-0.1,-0.20)   node [left] {$0$}; 
 \end{tikzpicture}
\end{center}
\caption{Space upper multifractal spectrum  of $\mu^\alpha$.}
\label{fig-stable}
\end{figure}
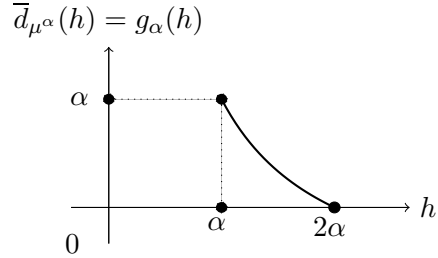
\end{center}

Our first result gives the possible values for the local dimensions of the occupation measure $\mu$ associated with a stable-like process $\cM$. 

\begin{theo}\label{theoexponent} Consider  a stable-like process $\cM$   associated to a non-decreasing mapping $\be$, as in Definition \ref{beta}, and the associated  occupation measure $\mu$. With probability $1$, for every $x\in \supp(\mu)$,  
$$\underline{\dim}(\mu,x) = \be(x)  \ \mbox{ and }  \ \overline{\dim}(\mu,x)\in [\be(x),2\be(x)].$$
\end{theo}
 Hence,  the support of the lower spectrum $\un d_\mu $ is random, depending on the trajectory of $\cM$.

The {\em space} lower multifractal spectrum is then quite easy to understand,  since the level set $\un E_\mu(\cO,h)$ contains either one point or is empty, depending on whether $h$ belongs to the closure of the range of the index process $\{\be(\cM_{t}) :  \cM_t\in \cO\}$ or not.    
Theorem \ref{theoexponent} indicates that the spectrum related to the upper local dimension  $\ov \dim(\mu,\cdot)$ should be more interesting. This is indeed the case, as resumed in Theorem \ref{theospace}.   Set 
\begin{align*}
\hat{g}_\al(h) :=\begin{cases}  \alpha\left(\frac{2\alpha}{h}-1\right) &\ \mbox{ when }h\in [\alpha,2\alpha ) , \\ -\infty & \mbox{ otherwise.} \end{cases}
\end{align*}
Note that the only difference between $g$ and $\hat{g}$ is at the value $h= 2\al$. 

\begin{definition}
 For every monotone c\`adl\`ag function $\Upsilon: \R^+\to \R$, we denote by $S(\Upsilon)$ the set of jumps of   $\Upsilon$.
\end{definition}

\mk
\begin{theo} \label{theospace} 
Set  the (at most countable) sets of real numbers  
 \begin{align}
 \label{defE}
 \mathcal{E}_1 &= \{\be(\cM_t) : t \in S(\cM) \mbox{ and } \be(\cM_t)\ge 2\be(\cM_{t-})  \}, \nonumber \\
  \cE_2&= \{2\be(\cM_{t-}) : t \in S(\cM) \mbox{ and } \be(\cM_t)\ge 2\be(\cM_{t-})   \},  \nonumber\\
 \cE& = \cE_1\cup \cE_2.
 \end{align}
 With probability 1, for  every non-trivial open interval $\cO\subset\rr$,  one has 
\begin{align}
\label{spacespectrum}
\un d_\mu(\cO,h)  & =     \begin{cases}  \ \ 0 &\mbox{ if } \, h \in\overline{\{\be(\cM_t): \cM_t\in \cO\}},\\
-\infty &\mbox{ otherwise,}
\end{cases}\end{align}
and for every $h \in \RR^+ \setminus \mathcal{E}$,
\begin{align}
 \label{spacespectrum2}
  \ov d_{\mu}(\cO,h) &=     \sup \, \big \{\hat{g}_\alpha(h) : \alpha\in  \{\be(\cM_t): \cM_t\in \cO\}  \big \}  .
\end{align} 

\end{theo}  
\begin{remark} If the range of $\be(\cdot)$ is included in, for example,  $[1/2, 9/10]$,   then the set of exceptional values $\cE  = \emptyset$ a.s. 
\end{remark}

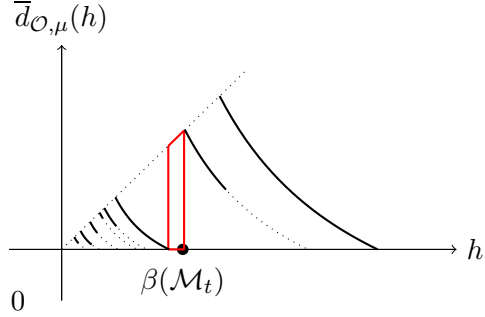
\begin{figure}
\begin{center}
\begin{tikzpicture}[xscale=3.5, yscale=3.4]
\draw [fill,->] (-0.2,-0)--(1.5,0) node [right] {$h$  };
\draw [fill,->] (0,-0.2)--(0,0.8) node [above] {$\ov d_{\cO, \mu }(h)$  };
   \draw [dotted, domain=0:0.7, color=black]  plot ({\x},  {\x});
\draw [fill] (-0.1,-0.20)   node [left] {$0$}; 

\draw  [fill] (0.46,0)circle [radius=0.02]    [fill]  (0.46,-0.03) node [below] {$\beta(\cM_t)$  };

\begin{scope}[xshift=-0.0cm,scale=0.75] 
\draw [fill, color = red, thick] (0.54,0.54)--(0.54,0) (0.54,0)--(0.62,0)  (0.62,0)--(0.62,0.62)   (0.54,0.54)--(0.62,0.62)    ; 
\end{scope}

 \draw [thick, domain=0.6:1.2, color=black]  plot ({\x},  {0.6*(2*0.6/\x-1)});

\begin{scope}[xshift=-0cm,scale=0.78] 
 \draw [ dotted,domain=0.6:1.2, color=black]  plot ({\x},  {0.6*(2*0.6/\x-1)});
 \draw [ thick,domain=0.6:0.8, color=black]  plot ({\x},  {0.6*(2*0.6/\x-1)});

\end{scope}

\begin{scope}[xshift=-0.0cm,scale=0.08] 
 \draw [ dotted,domain=0.6:1.2, color=black]  plot ({\x},  {0.6*(2*0.6/\x-1)});
 \draw [ thick,domain=0.6:0.8, color=black]  plot ({\x},  {0.6*(2*0.6/\x-1)});
\end{scope}

\begin{scope}[xshift=-0.0cm,scale=0.12] 
 \draw [ dotted,domain=0.6:1.2, color=black]  plot ({\x},  {0.6*(2*0.6/\x-1)});
 \draw [ thick,domain=0.6:0.95, color=black]  plot ({\x},  {0.6*(2*0.6/\x-1)});
\end{scope}

\begin{scope}[xshift=-0.0cm,scale=0.18] 
 \draw [ dotted,domain=0.6:1.2, color=black]  plot ({\x},  {0.6*(2*0.6/\x-1)});
 \draw [ thick,domain=0.6:0.75, color=black]  plot ({\x},  {0.6*(2*0.6/\x-1)});
\end{scope}

\begin{scope}[xshift=-0.0cm,scale=0.24] 
 \draw [ dotted,domain=0.6:1.2, color=black]  plot ({\x},  {0.6*(2*0.6/\x-1)});
 \draw [ thick,domain=0.6:0.65, color=black]  plot ({\x},  {0.6*(2*0.6/\x-1)});
\end{scope}

\begin{scope}[xshift=-0.0cm,scale=0.27] 
 \draw [ dotted,domain=0.6:1.2, color=black]  plot ({\x},  {0.6*(2*0.6/\x-1)});
 \draw [ thick,domain=0.6:0.77, color=black]  plot ({\x},  {0.6*(2*0.6/\x-1)});
\end{scope}

\begin{scope}[xshift=-0.0cm,scale=0.34] 
 \draw [dotted,domain=0.6:1.2, color=black]  plot ({\x},  {0.6*(2*0.6/\x-1)});
 \draw [ thick,domain=0.6:1.2, color=black]  plot ({\x},  {0.6*(2*0.6/\x-1)});
\end{scope}
 
\end{tikzpicture}
\end{center}
\caption{Space upper multifractal spectrum of $\mu $ on an open set $\cO$. The spectrum (in thick) is obtained as the supremum of a random  countable number of functions of the form $g_\alpha(h)$, for the values $\alpha \in \beta(\cM_{t-})$, $\cM_t\in \mathcal{O}$.   It may happen that there is a hole in the support of  $\ov d_{\cO, \mu }$ (in red in the figure). In this case, the value of $\ov d_{\cO, \mu }$ at $\beta(\cM_t)$ is either 0 or $-\infty$. } 
\label{fig-markov1} 
\end{figure} 


First, one shall notice that both spectra are random, depending on the trajectory {\em and} on the interval $\cO$. In this sense, $\un d_\mu(\cO,\cdot) $ and $\ov d_\mu(\cO,\cdot) $  are inhomogeneous, contrarily to what happens for the occupation measure $\mu^\alpha$ of $\alpha$-stable subordinators (the spectra do not depend on $\cO$).

One shall interpret the space upper spectrum as the supremum of an infinite number of space multifractal spectra of "locally $\alpha$-stable processes" for all values $\alpha\in  \{\be(\cM_t): \cM_t\in \cO\}$. This formula finds its origin in the fact that locally, $\cM$ behaves around each continuous time $t$ as an $\alpha$-stable process with $\alpha = \be(\cM_t)$.

The allure of a typical space upper multifractal spectrum is depicted in   Figure \ref{fig-markov1}. This  shape is very unusual in the literature. 

First, observe that, since $\beta$ and $\cM$ are increasing maps,  when $t_0\in S(\cM)$ is a jump time for $\cM$, then the "local" index of $\cM$ jumps at $t_0$  from $\be(\cM_{t_0-})$ to $\be(\cM_{t_0})$, and for $t\geq t_0$, the only possibility to have $ \ov d_{\mu}(\cO,\cM_t) = \be(\cM_{t_0})$ is when $t=t_0$. Similarly, when $t<t_0$, it is not possible to have $ \ov d_{\mu}(\cO,\cM_t) = 2\be(\cM_{t_0-})$.

In particular    there may be a "hole" in the support of $\ov d_\mu(\cO,\cdot)$. Indeed, a quick analysis of the functions $g_\alpha(\cdot)$ shows that this happens when there is a   time $t_0$ such that  $\be(\cM_t) > 2 \be(\cM_{t-})$, which occurs with positive probability for functions $\be(\cdot)$ satisfying $2\e_0< 1-\e_0$.  

All this explains the set of exceptional points $\mathcal{E}$ in Theorem \ref{theospace}. We deal with these exceptional points in the following theorem, whose statement is rather long but whose proof    follows directly from a careful analysis of the previous results.

\begin{theo} \label{theospaceextreme} 
With probability 1, for  every non-trivial open interval $\cO\subset\rr$,  when $h \in \mathcal{E}$, three cases may occur. 
\begin{enumerate}
\item  $h=\be(\cM_t)>2\be(\cM_{t-})$.  If $\cM_t\notin\cO$,  $\ov d_\mu(\cO, h)= -\infty$, otherwise one has 
$$ \ov d_{\mu}(\cO,h)= \begin{cases} 
 \ \ \ \ \ \ \ 0  & \mbox{ if }  \ov \dim(\mu,\cM_t) = h,\\  
 \ \ \ \ -\infty &  \mbox{ if }   \ov \dim(\mu,\cM_t) > h.
 \end{cases}$$ 
 
 \sk
 
 \item  $h=2\be(\cM_{t-})<\be(\cM_t)$. If $\cM_{t-}\notin\cO$,  $\ov d_\mu(\cO, h)= -\infty$, otherwise one has 
$$ \ov d_{\mu}(\cO,h)=\begin{cases} \ \ \ \ \ \ \ \ 0  & \mbox{ if }   \ov\dim(\mu, \cM_{t-}) = h,  \\ 
 \ \ \ \ \ \ \ -\infty  & \mbox{ if } \ov\dim (\mu, \cM_{t-}) <  h.\end{cases}$$

\item $h = \be(\cM_{t}) = 2\be(\cM_{t-})$. If $\{\cM_t,\cM_{t-}\}\in \cO$,  one has
 $$ \ov d_{\mu}(\cO,h)=\begin{cases}  \ \ \ \ \ \ 0  & \mbox{ if }   \ov\dim(\mu, \cM_{t-}) = h  \mbox{ or }   \ov\dim(\mu, \cM_{t}) = h, \\ 
 \ \ \ -\infty  & \mbox{ if }   \ov\dim(\mu, \cM_{t-}) < h  \mbox{ and  }   \ov\dim(\mu, \cM_{t}) > h.\end{cases}$$ 
If only one of $\cM_t$ and $\cM_{t-}$ belongs to $\cO$ (say, $\cM_{t-}$), one has 
    $$ \ov d_{\mu}(\cO,h)=\begin{cases}  \ \ \ 0  & \mbox{ if }   \ov\dim(\mu, \cM_{t-}) = h,  \\ 
 -\infty  & \mbox{ if }   \ov\dim(\mu, \cM_{t-}) < h.\end{cases}$$ 
 If neither $\cM_t$ nor $\cM_{t-}$ belongs to $\cO$,  one has $\ov d_{\mu}(\cO,h)= -\infty$.
\end{enumerate}

\end{theo}  

 \medskip

It is also interesting to consider the {\em time multifractal spectra} related to $\cL^\alpha$ and $\cM$, which describes the variation along time of the Hausdorff dimension of the set of times $t$ such that $\mu^\alpha$ (or $\mu$) has a local dimension $h$ at $x =  \cL^\alpha_t$ (or $\cM_t$). By abuse of language, one says in this case  that $\mu$ (or $\mu^\alpha$) has a dimension $h$ at $t$. For this let us introduce other level sets.

\begin{definition}\label{defi4} 
  For every open set $\cO\subset [0,1]$, set
\begin{align*}
  \ov E_{\mu^\alpha}^t(\cO, h)&=\{t\in \cO : \ov \dim(\mu^\alpha,\cL^\alpha_t) =h \},  \\
   \ov E_{\mu}^t(\cO, h)&=\{t\in \cO : \ov \dim(\mu^\alpha,\cM_t) =h \}, 
\end{align*}
and the similar quantities for lower local dimensions $  \un E^t_{\mu^\alpha}(\cO, h)$ and $ \un E_{\mu}^t(\cO, h)$.
 The corresponding {\em time multifractal spectra}  of $\mu^\alpha$ and $\mu $ are
\begin{align*}
 \overline{d}_{\mu^\alpha}^t(\cO,\cdot):  &  \ h\mapsto \dim_{\cH} \ov E^t_{\mu^\alpha}(\cO, h),\\  
 \overline{d}_{\mu }^t(\cO,\cdot):  & \ h\mapsto \dim_{\cH} \ov E^t_{\mu }(\cO, h),
  \end{align*}
and the similar quantities for lower local dimensions.
\end{definition}

In the case of a stable subordinator, Hu and Taylor prove that a.s., 
\begin{equation}
\label{spec-stable3}
\overline{d}_{\mu^\alpha}^t(\cO,h) =  \frac{g_\alpha(h)}{\alpha} =\begin{cases}  2\alpha/h-1 & \mbox{ when }h\in [\alpha,2\alpha],\\ -\infty & \mbox{ otherwise}.\end{cases}
\end{equation}
 This time upper multifractal spectrum is {\em homogeneous}, in the sense that {\em it does not depend on the choice of $\cO$}. In this article, we also compute the time   multifractal spectra of $\mu$.

\begin{theo} \label{theotime} 
Set the (at most countable) set of real numbers $$\cE' = \{ \be(\cM_t) : t\in S(\cM) \mbox{ and }  \be(\cM_t) \ge 2\be(\cM_{t-})  \}. $$
With probability 1, for  every non-trivial open interval $\cO\subset[0,1]$,   
\begin{align}
\label{timespectrum}
\un d^t_\mu(\cO,h)   =     \begin{cases}  \ \ 0 &\mbox{ if } \, h \in {\{\be(\cM_t): t\in \cO\}},\\
-\infty &\mbox{ otherwise,}
\end{cases}
\end{align}
and for every $h \in \cE'$, 
\begin{align}
\label{timespectrum2}
 \ov d^t_{\mu}(\cO,h) =     \sup \, \left \{ \frac{\hat{g}_\alpha(h)}{\alpha}  : \alpha\in  \{\be(\cM_t): t\in \cO\}  \right \}  .
\end{align}
\end{theo}  
\begin{remark}
Note that $\cE'\subset \cE$.  If, for example,  the range of $\be(\cdot)$ is included in $[1/2, 9/10]$, then almost surely $\cE'= \emptyset$.
\end{remark}
The first part is trivial. Observe that there is a subtle difference between \eqref{timespectrum}  and \eqref{spacespectrum}, since at each jump time $t$ for $\cM$, there is no $s\in \R$ such that $\cM_s=\cM_{t-}$.  

As for the space multifractal spectrum, there is some uncertainty about the value of $ \ov d^t_{\mu}(\cO,\cM_t)$ at the jump  times $   t \in S(\cM)  $, which is dealt with in the following theorem.

\begin{theo} \label{theospaceextreme2} 
With probability 1, for  every non-trivial open interval $\cO\subset\rr$,  for every $h \in  \{\be(\cM_t) : t \in S(\cM) \} $, writing $h=\beta(\cM_t)$ , one has
$$ \ov d^t_{\mu}(\cO,h)=   \frac{\ov d_{\mu}(\cO,h) }{\be(\cM_{t-})} .$$ 
\end{theo}  

\begin{figure}
\begin{center}
\begin{tikzpicture}[xscale=3.5, yscale=3.4]
\draw [fill,->] (-0.2,-0)--(1.5,0) node [right] {$h$  };
\draw [fill,->] (0,-0.2)--(0,1.2) node [above] {$\ov d^t_{\mu }(\cO, h)$  };

\draw [fill] (-0.1,-0.20)   node [left] {$0$}; 

 \draw [ dotted,domain=0.24:0.48, color=black]  plot ({\x},  {(2*0.24/\x-1)});
 \draw [ thick,domain=0.24:0.28, color=black]  plot ({\x},  {(2*0.24/\x-1)});
 
 \draw [ dotted,domain=0.28:0.56, color=black]  plot ({\x},  {(2*0.28/\x-1)});
 \draw [ thick,domain=0.28:0.3, color=black]  plot ({\x},  { (2*0.28/\x-1)});

 \draw [ dotted,domain=0.3:0.6, color=black]  plot ({\x},  { (2*0.3/\x-1)});
 \draw [ thick,domain=0.3:0.4, color=black]  plot ({\x},  { (2*0.3/\x-1)});

 \draw [ dotted,domain=0.4:0.8, color=black]  plot ({\x},  {(2*0.4/\x-1)});
 \draw [ thick,domain=0.4:0.55, color=black]  plot ({\x},  { (2*0.4/\x-1)});
 
 \draw [dotted,domain=0.55:1.1, color=black]  plot ({\x},  { (2*0.55/\x-1)});
 \draw [ thick,domain=0.55:0.62 , color=black]  plot ({\x},  { (2*0.55/\x-1)});

 \draw [ dotted,domain=0.62:1.24, color=black]  plot ({\x},  { (2*0.62/\x-1)});
 \draw [ thick,domain=0.62:0.7, color=black]  plot ({\x},  { 2*0.62/\x-1)});

 \draw [ dotted,domain=0.7:1.4, color=black]  plot ({\x},  { (2*0.7/\x-1)});
 \draw [ thick,domain=0.7:1.4, color=black]  plot ({\x},  { (2*0.7/\x-1)}); 
 
\end{tikzpicture}
\end{center}
\caption{Time upper multifractal spectrum of $\mu $ on an open set $\cO$. The spectrum (in thick) is obtained as the supremum of a random  countable number of functions of the form $\frac{g_\alpha(h)}{\alpha}$ (drawn using dotted graphs), for the values $\alpha \in \beta(\cM_{t})$, $t\in \mathcal{O}$.    } 
\label{fig-markov2} 
\end{figure}
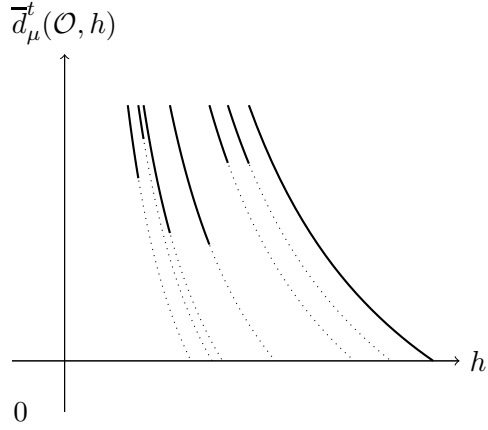

The correspondance between \eqref{spec-stable2} and  \eqref{spec-stable3} follows from the fact that for  every $\al$-stable subordinator,  almost surely for each measurable set $E\subset[0,1]$, 
\begin{align}\label{uniformrange}
\dimh(\cL^\alpha(E)) = \al\cdot\dimh(E).
\end{align}
Up to a countable number of points,  writing $\cO=(g,d)$, one has the equality  $\cL^\alpha \Big(\ov E_{\mu^\alpha}^t(\cO, h)\Big) = \ov E_{\mu^\alpha}\left( (\cL^\alpha_g, \cL^\al_d),h\right)$. The method developed by Hu and  Taylor consists first in proving  \eqref{spec-stable3}, and then in applying  \eqref{uniformrange} to get   \eqref{spec-stable2}.  

\smallskip

Following these lines, we   start by proving  Theorem \ref{theotime}. The original methods by Hu and Taylor do not  extend here, and   an alternative way to compute the time multifractal spectrum of $\mu$ is needed. For this,    some scenario leading to the fact that $\mu$ has exactly an upper local dimension equal to $h$ at $x=\cM_t$ is identified. More precisely, it will be proved that $d_\mu(\cO,\cM_t)=h$ when $t$ is infinitely many times very closely surrounded by two "large" jump times for the Poisson point process involved in the construction of $\cM$. Using this property, we build in Section \ref{seclowerbound} a (random) Cantor set of such times $t$ with the suitable Hausdorff dimension. The difficulty lies in the fact that the expected Hausdorff dimension is random and depends on the interval we are working on.

\smallskip

Further, it is natural to look for a uniform dimension formula such as \eqref{uniformrange} for the stable-like jump diffusion process $\cM$.  As is pointed out by Hu and Taylor \cite{hu2000occupation} (see also page 94 of \cite{bertoin1996book}),  as long as the Laplace exponent of a general subordinator oscillates at infinity, i.e. it exhibits different power laws at infinity,  one can never expect  such an identity to hold.  Nevertheless, using an   argument  based on coupling and   time change, we  find   upper and lower dimension bounds for  images of sets by  stable-like processes.
\begin{theo}
\label{th_dim_images}
Almost surely, for every mesurable set $E\subset [0,1]$, one has
\begin{equation}
\label{eqfinal}
 \dimh \cM(E)  \in \dimh(E)  \cdot  \Big[\inf_{t\in E} \beta(\cM(t)),\sup_{t\in E} \beta(\cM(t-)) \Big].
 \end{equation}
 
Moreover, if  the set $E$ satisfies that for every non-trivial subinterval $\cO\subset [0,1]$, $\dimh (E) = \dimh (E\cap \cO)$, then
$$ \dimh \cM(E)  = \dimh(E)  \cdot  \sup_{t\in E} \beta(\cM(t-))  .$$
\end{theo}
These results  are fine enough for us to deduce Theorem \ref{theospace} from Theorem \ref{theotime}, even when the stability parameter of these processes varies as time passes. Also, this theorem solves partially a question left open in \cite{yang2015range}.

\sk

Let us end this introduction with a proposition describing the typical behavior of the occupation measure $\mu$.   We skip the proof since it can be deduced by  adapting the arguments in the proof of our Theorem  \ref{theogeneral}.
\begin{pro}[Typical Behavior] With probability $1$,  one has:
\begin{itemize}
\item for Lebesgue-almost every time $t\in[0,1]$,   $  \dim(\mu,\cM_t) = \be(\cM_t)$.  

\item for $\mu$-almost every point $x\in \R$, $ \dim(\mu,x)=\be(x)$.
\end{itemize}
\end{pro}

%

\sk
The rest of this paper is organized as follows.  We start by recalling  basic properties of the stable-like processes  in Section \ref{secprelim}. The local dimensions of the occupation measure (Theorem \ref{theoexponent}) are studied in Section \ref{sechold}. The time spectrum (Theorems   \ref{theotime} and \ref{theospaceextreme2}) is obtained in Section \ref{sectimegeneral} using a general result  (Theorem \ref{theogeneral}), whose proof is given in Sections \ref{secupperbound} and \ref{seclowerbound}.  Finally, the space spectrum (Theorem \ref{theospace}) is dealt with  in Section \ref{sec_space}, together with  the dimension of images of arbitrary sets by stable-like processes (Theorem \ref{th_dim_images}).

\section{Preliminaries}\label{secprelim}

First of all, stable-like processes admit a Poisson representation which was regularly used to study path properties of such processes, see for instance \cite{barral2010increasing,yang2015jumpdiffusion,yang2015range}. Let us recall this representation and a coupling associated with it which will be useful for our purposes.

\sk
Let $N(dt,dz)$ be a Poisson measure on $\R^+\times \R$  with intensity $dt\otimes dz/z^2$. Such a measure can be constructed from a Poisson point process which is  the set of jumps of a L\'evy process with triplet $(0,0,dz/z^2)$, see for instance Chapter 2 of \cite{applebaum2009}. We denote $\cF_t = \sigma \( \{ N(A): A\in  \cB([0,t]\times[0,+\infty)) \} \) $. 

Recall the definition of a stable-like process and formula \eqref{defstablelike}.  The existence and uniqueness of such jump diffusion processes is classical and recalled in the next proposition. Observe that by the substitution $u=z^{1/\be(x)}$(for each fixed $x$), the generator of a stable-like process  is rewritten as
\begin{align}\label{generator}
\mathfrak{L}f(x)= \int_0^1 \( f(x+z^{\frac{1}{\be(x)}})-f(x)\) \frac{dz}{z^2}.
\end{align}

\begin{pro}\label{coupling} Let $N$ be as in the last paragraph.
\begin{enumerate}
\item There exists a unique c\`adl\`ag $(\cF_t)_{t\in[0,1]}$-adapted solution to
\begin{equation}
\label{defM}
\cM_t = \int_0^t\int_0^1 z^{\frac{1}{\be(\cM_{u-})}} N(du,dz).
\end{equation}
Furthermore, $\cM$ is an increasing strong Markov process with generator $\mathfrak{L}$ given by  \eqref{generator}. 
\item For every $\al\in(0,1)$, we define 
\begin{align}\label{defLalpha}
\cL^\al_t = \int_0^t\int_0^1 z^{\frac{1}{\al}} N(du,dz).
\end{align}
Then for all $\al\in(0,1)$, $\{\cL^\al_t, t\in[0,1]\}$ is an $\al$-stable subordinator whose jumps larger than $1$ are truncated.
\end{enumerate}
\end{pro}

Classical arguments based on Gronwall inequality and Picard iteration yield the first item. For a proof, see Proposition 13 of \cite{barral2010increasing} or Proposition 2.1-2.3 of \cite{yang2015jumpdiffusion} with some slight modifications. The second item is standard, see for instance Section 2.3 of \cite{applebaum2009}.

\begin{remark}
\label{remcoupling}
Recall that   $S(\Upsilon)$  is the set of jumps of    a  monotone c\`adl\`ag function $\Upsilon: [0,1]\to \R^+$.

Observe that by construction, almost surely, the processes $\cM$ and the family of L\'evy processes $(\cL^\alpha)_{\al\in(0,1)}$ are purely discontinuous,  increasing, with finite variation, and that they jump at the same times, i.e. $S(\cM)= S(\cL^\alpha)$.
\end{remark}

Next observation is key for the study of the local dimensions of $\mu$.

\begin{pro}\label{procompare} Consider the process $\cM$ and $\cL^\al$ for all $\al\in(0,1)$ introduced in Proposition \ref{coupling}. Almost surely, for all $0\le s\le t\le 1$, 
\begin{align*}
0\le \cL^{\be(\cM_s)}_t - \cL^{\be(\cM_s)}_s \le \cM_t - \cM_s \le \cL^{\be(\cM_{t-})}_t - \cL^{\be(\cM_{t-})}_s.
\end{align*}
\end{pro}
This is intuitively true because we construct simultaneously $\cM$ and $\cL^\al$ such that they jump at the same times, and the jump size of $\cL^{\al}$ is always larger than $\cL^{\al'}$ whenever $\al>\al'$. See Proposition 14 of \cite{barral2010increasing} for a proof.

\section{Local dimensions of $\mu$ : Proof of Theorem \ref{theoexponent}}\label{sechold}
Observe that almost surely, for all $\al\in\Q\cap(0,1)$,  formula \eqref{spec-stable} is true. This, together with Proposition \ref{procompare}, leads to the local dimension of $\mu$.

\sk
\begin{preuve}{\it \ of Theorem \ref{theoexponent} : }
Three cases may occur.
\begin{enumerate}
\item $x=\cM_t$ with $t$ a continuous time of $\cM$.  
Due to the coupling in their construction (Proposition \ref{coupling}), almost surely,  for every $\alpha$, the process $ \cL^\alpha$ is also continuous at $t$. 

By continuity, for arbitrary rational numbers $\al,\al'\in(0,1)$ satisfying $\al<\be(\cM_t)<\al'$, there exists a small $\de>0$ such that  for all $s\in(t-\de,t+\de)$, $\al<\be(\cM_s)<\al'$. Using    the occupation measure $\mu^\al$ of the process $\cL^\al$, and applying Proposition \ref{procompare} to $\cL^\al$ and $\cL^{\al'}$, one gets  when  $r $ is small
\begin{align}\label{measurecompare}
\mu^{\al'}  \big( (\cL^{\al'}_t -r, \cL^{\al'}_t+r)  \big)\le \mu \big((\cM_t-r,\cM_t+r)  \big)\le \mu^\al \big((\cL^\al_t-r,\cL^\al_t+r) \big)
\end{align}

By formula \eqref{spec-stable} for the lower and upper local dimensions of $\mu^\alpha$, for all small $\e>0$, almost surely, one has for   $r$ small enough that 
\begin{align*}
\alpha-\e \le  \frac{\log \mu^\al \big ((\cL^\al_t-r,\cL^\al_t+r)\big)}{\log(r)} \le 2 \al+\e,
\end{align*}
and the same for $\alpha'$. Hence \begin{align*}
\frac{\log \mu \big ((\cM_t-r,\cM_t+r) \big )}{\log(r)} \ge \frac{\log\mu^\al \big ((\cL^\al_t-r,\cL^\al_t+r) \big )}{\log(r)} \ge \al-\e.
\end{align*}
and \begin{align*}
\frac{\log \mu\big ((\cM_t-r,\cM_t+r)\big)}{\log(r)} \le \frac{\log \mu^{\al'}\big ((\cL^{\al'}_t-r,\cL^{\al'}_t+r)\big)}{\log(r)} \le 2\al'+\e.
\end{align*}
 Therefore, $ \alpha-\e \leq \underline{\dim}(\mu,x) \leq \ov \dim(\mu,x)\le 2 \al'+\e$.

On the other hand, still by formula \eqref{spec-stable}, $\un\dim(\mu^{\al'},x) =\alpha'$, so there exists a sequence $(r_n)$ converging to 0 such that 
\begin{align*}
\alpha'-\e \le  \frac{\log \mu^\al \big ((\cL^{\al'}_t-r_n,\cL^{\al'}_t+r_n)\big)}{\log(r_n)} \le \al'+\e,
\end{align*}
so $ \alpha-\e \leq \underline{\dim}(\mu,x) \leq \alpha'+\ep$.

Letting $\ep$ tend to zero and $\alpha$, $\alpha'$ tend to $\beta(\cM_t)$ with rational values yields $ \beta(\cM_t)=\underline{\dim}(\mu,x) \leq \ov \dim(\mu,x)\le 2\beta(\cM_t)$.

\mk

\item $x=\cM_t$ with  $t$ a jump time for $\cM$.   Observe that in this case $\mu\big((x-r,x+r) \big)= \mu\big((\cM_t, \cM_t+r)\big)$ for $r>0$ small enough. For arbitrary rational numbers $\al<\be(\cM_t)<\al'$, the inequality \eqref{measurecompare} is straightforward using Proposition \ref{procompare}. We follow the same lines in the first case to obtain  the desired result.
\mk

\item  $x=\cM_{t-}$ with $t$ a jump time for $\cM$.    Now, $\mu\big((x-r,x+r) \big)= \mu\big((\cM_{t-}-r, \cM_{t-}) \big)$ for $r>0$ small. Then the  proof goes like the previous items. \end{enumerate}
\end{preuve}

Let us end this Section with the proof of the easier part of Theorem \ref{theospace} :  space lower multifractal spectrum of $\mu$. 
\sk

\begin{preuve}{\it \ of  Formula \eqref{spacespectrum} of Theorem \ref{theospace} :}
As noticed in Remark \ref{remcoupling}, $t\mapsto\be(\cM_t)$ is increasing due to   the monotonicity of  $\cM$ and $\be$.  Hence each level set of $t\mapsto\be(\cM_t)$ contains at most one point. This means that for each open interval $\cO$ that intersects $\supp(\mu)$, 
\begin{align*}
\un E_\mu(\cO, h) & = \{x\in\supp(\mu)\cap\cO : \be(x)=h \}\\ 
&  = \begin{cases} \{\be(\cM_t)\} &\mbox{ if } h= \be(\cM_t) \mbox{ for some } t \mbox{ with } \cM_t\in\cO,\\
\{\be(\cM_{t-})\} &\mbox{ if } h= \be(\cM_{t-}) \mbox{ for some } t  \mbox{ with } \cM_t\in\cO, \\
\emptyset &\mbox{ if } h\not\in \overline{\{\be(\cM_t): \cM_t\in\cO\}},
\end{cases}
\end{align*}
which completes the proof.
\end{preuve}

\section{A   general result to get the  time spectrum  (Theorems \ref{theotime} and \ref{theospaceextreme2})} \label{sectimegeneral}

Let us present a general result, proved in Sections \ref{secupperbound} and \ref{seclowerbound}.  This theorem gives the dimension of the random set of times $t$ where the local dimension mapping  $s\mapsto \ov \dim (\mu, \cM_s)$ coincides with a given function.  The remarkable feature of this theorem is that it allows to determine these dimensions for all the monotone  c\`adl\`ag function simultaneously,  with probability one.

\begin{theo}
\label{theogeneral}
For every non-increasing c\`adl\`ag function $\Upsilon: [0,1]\to [1,2]$ and every open interval  $\cO\subset [0,1]$, let us define  $\Upsilon_{\min}=\inf_{u\in\cO }\Upsilon(u)$ and the sets 
\begin{align*}
\ov E_\mu ^t(\cO,\Upsilon) &= \llb t\in \cO : \ov \dim(\mu,\cM_t) = \Upsilon(t)\be(\cM_t) \rrb, \\
\ov E_\mu ^{t,\geq}(\cO,\Upsilon)  &= \llb t\in \cO\  : \ov \dim(\mu,\cM_t) \ge \Upsilon(t)\be(\cM_t) \rrb.
\end{align*}

With probability $1$, for every non-increasing c\`adl\`ag function $\Upsilon: [0,1]\to [1,2]$ and every open interval  $\cO\subset [0,1]$,  we have
\begin{equation}
\label{eqth7}
\dimh \ov E_\mu ^t(\cO,\Upsilon)   = \dimh \ov E_\mu ^{t,\geq}(\cO,\Upsilon) = \dimh \ov E_\mu ^{t,\geq}(\cO,\Upsilon_{\min}) = \frac{2}{\Upsilon_{\min}} -1.
\end{equation}
\end{theo}
%

\noindent The notation $ \ov E_\mu ^{t,\geq}(\cO,\Upsilon_{\min})$ means that we consider the constant function $\Upsilon \equiv \Upsilon_{\min}$.

\subsection{Proof for the time upper multifractal spectrum}

Let us  explain why Theorems \ref{theoexponent} and   \ref{theogeneral} together imply Theorems \ref{theotime} and \ref{theospaceextreme2}. One wants to prove formula  \eqref{timespectrum2}, which can be rewritten as
\begin{equation}
\label{timespectrum22}
 \ov d^t_{\mu}(\cO,h) =     \sup \, \left \{ \frac{2\alpha}{h}-1\in[0,1)  : \alpha\in  \{\be(\cM_t): t\in \cO\}  \right \}  .
 \end{equation}
 
 \mk
 
  One combines Theorem \ref{theoexponent} and Theorem \ref{theogeneral} with the family of functions $\{\Upsilon^h(t) = h/\be(\cM_t) : h\ge 0 \}$.   With probability one, these functions are all  c\`adl\`ag decreasing. Observe that for every open interval $\cO$, 
\begin{align*}
\dimh \ov E_\mu ^t(\cO,\Upsilon^h)   & =  \dimh  \llb t\in \cO : \ov \dim(\mu,\cM_t) = \Upsilon^h(t)\be(\cM_t) \rrb \\
&=    \dimh  \llb t\in \cO : \ov \dim(\mu,\cM_t) = h \rrb  \\
&=  \ov d^t_\mu(\cO,h).
  \end{align*}

  We prove now that formula \eqref{eqth7} applied to the family $\{\Upsilon^h: h\ge 0 \}$  implies formula \eqref{timespectrum22}.
  Several cases  may occur according to the value of $h$. 
  
  \medskip
  
$\bullet$ {\bf If $\{1\}  \not\subset \Upsilon^h(\cO)$:}
           \begin{enumerate}
           \item {\bf First case:} For all $t\in\cO$, $\Upsilon^h(t)>1.   $
                    \begin{itemize}
  \smallskip                    \item If  $\inf_{t\in\cO} \Upsilon^h(t)\ge 2$,  then for all $t\in\cO$, $\Upsilon^h(t)>2$.  Theorem \ref{theoexponent} entails  $\ov E_\mu ^t(\cO,\Upsilon^h)  = \emptyset$.   So  $\ov d^t_\mu(\cO,h)=\dimh \ov E_\mu ^t(\cO,\Upsilon^h)   =-\infty$, which coincides with   \eqref{timespectrum22}.
  
    \smallskip                    \item If   $\inf_{t\in\cO} \Upsilon^h(t) < 2$,  consider the entrance time in $(1,2)$ by $\Upsilon^h$  $$\tau  = \inf\{ t\in\cO : \Upsilon^h(t)<2 \}. $$
  By construction, $\forall\, t\in(\tau,\infty)\cap \cO$,  $\Upsilon^{h}(t)\in (1,2)$. By Theorems \ref{theoexponent} and \ref{theogeneral},  one gets
                    \begin{align*} 
                    \ov d^t_\mu(\cO,h) &= \ov d^t_\mu((\tau,\infty)\cap\cO ,h) = \frac{2}{\inf_{t\in(\tau,\infty)\cap\cO} \Upsilon^h(t)} -1 \\
                    & = \frac{2\sup_{t\in(\tau,\infty)\cap\cO}\be(\cM_t)}{h} -1,
                    \end{align*}
                    which coincides with \eqref{timespectrum22}. 
                    \end{itemize}

   \smallskip  \smallskip          \item  {\bf Second case:} There exists $t\in\cO$, such that $\Upsilon^h(t)<1$. Define the passage  time of $(-\infty, 1)$ by $\Upsilon^h$  as
   $$\sigma = \inf\{t\in\cO : \Upsilon^h(t)<1 \}. $$  
   
           \begin{itemize}
    \smallskip           \item If $\sigma$ is the left endpoint of $\cO$, then for all $t\in\cO$, $\Upsilon^h(t)<1$.  Theorem \ref{theoexponent} yields $\ov E_\mu ^t(\cO,\Upsilon^h)   = \emptyset$.   Again, this gives  $\ov d^t_\mu(\cO,h)= -\infty$, which coincides with   \eqref{timespectrum22}.

    \smallskip           \item If $\sigma$ belongs to the open interval $\cO$,  the proof goes along the same lines as in item 1. replacing  $\cO$  by $(-\infty,\sigma)\cap \cO$. 
           \end{itemize}
                    
           \end{enumerate}
         
         \medskip  
$\bullet$ {\bf If $\{1\}  \subset \Upsilon^h(\cO)$:}
Let $t_0\in\cO$ be the unique time such that  $\Upsilon^h(t_0) =1$, i.e. $h = \be(\cM_{t_0})$.   One distinguishes different cases according to the behavior of  $t\mapsto \ \cM_t$ at $t_0$. 
  \begin{enumerate}
  \smallskip  \item {\bf If $ \cM$ is continuous at $t_0$:}  $\be(\cM_{\cdot})$ is also continuous at $t_0$. By definition, the entrance time $\tau$ satisfies $\tau<t_0$.    Theorem \ref{theogeneral} entails
  \begin{align*}\ov d^t_\mu(\cO,h)  = \ov d^t_\mu((\tau,t_0)\cap\cO,h)  =  \frac{2}{\inf_{t\in(\tau,t_0)\cap\cO}\Upsilon^h(t)} -1 =1,
  \end{align*}  which coincides with \eqref{timespectrum22}. 

  \smallskip  \item  {\bf If $t_0$ is a jump time for $\cM$ and  $\be(\cM_{t_0-})< h = \be(\cM_{t_0}) < 2\be(\cM_{t_0-})$:}  Then, using that   $\Upsilon^h(t_0) = 1$, one deduces that  $0< \Upsilon^h(t_0-) - \Upsilon^h(t_0) <1$, which implies $\inf_{t\in (\tau, t_0)\cap\cO} \Upsilon^h(t) < 2$.  The same computation as in item 1.  with $\cO$ replaced by $(\tau,t_0)\cap\cO$ yields formula   \eqref{timespectrum22}. 
  
  \smallskip  \item  {\bf If $t_0$ is a jump time for $\cM$ and  $ h = \be(\cM_{t_0}) \geq 2\be(\cM_{t_0-})$:}  Then $\Upsilon^h(t_0-) \ge 2$, thus $\tau = t_0$.  One has 
  \begin{align*}
 \ov E^t(\cO, h) = \begin{cases}  \{t_0\} & \mbox{ if } \ov\dim(\mu, \cM_{t_0}) = h,  \\ \emptyset & \mbox{ otherwise. } \end{cases}
 \end{align*} 
This last formula coincides with the one claimed by Theorem \ref{theospaceextreme2}.
  \end{enumerate}

    \sk

\subsection{Reduction of the problem}

Observing that we have the obvious inclusion  $ \ov E_\mu ^t(\cO,\Upsilon)   \subset \ov E_\mu ^{t,\geq}(\cO,\Upsilon) \subset  \ov E_\mu ^{t,\geq}(\cO,\Upsilon_{\min})   $, we   proceed in two parts: 
\begin{itemize}
\item
first, in Section \ref{secup}, we show that 
\begin{equation}\label{reduc1}
\dimh\ov E_\mu ^{t,\geq}(\cO,\Upsilon_{\min}) \le \frac{2}{\Upsilon_{\min}} -1
\end{equation}  simultaneously for all $\Upsilon$ and $\cO$.

In order to get \eqref{reduc1}, it is equivalent to show that, almost surely, for each $\ga\in[1,2]$ and open interval $(a,b)\subset [0,1]$ with rational endpoints, 
   $$
\dimh \ov E_\mu ^{t,\geq}\big((a,b),\gamma\big) \le \frac{2}{\gamma}-1.
$$
   We will actually prove that for $\gamma\in(1,2)$,  almost surely,
   \begin{equation}\label{reduc2}
\dimh \ov E_\mu ^{t,\geq}\big((0,1),\gamma\big) \le \frac{2}{\gamma}-1.
\end{equation}
   The extension to  arbitrary $a,b\in[0,1]\cap\Q$ and  $\gamma\in\{1,2\}$ is straightforward.

\sk
\item
second,  in Section \ref{seclowerbound}, we complete the result by proving that 
\begin{equation}\label{reduclow1}
\dimh \ov E_\mu ^t(\cO,\Upsilon)  \ge \frac{2}{\Upsilon_{\min}} -1
\end{equation}  also simultaneously for all $\Upsilon$ and $\cO$, almost surely.
It is also enough to get the result for $\cO=(0,1)$. 

\end{itemize}

\section{Proof of Theorem \ref{theogeneral} : upper bound}\label{secupperbound} 
\label{secup}

Our aim is to prove \eqref{reduc2}. For notational simplicity, we write $\ov E_\mu ^{t,\geq}\big((0,1),\gamma\big) =    E(\gamma)$.  
Let us first observe that   the family of sets $\{E(\gamma), \gamma\in[1,2]\}$  is non-increasing  with respect to $\gamma$.  Recall that $\ga\in(1,2)$ throughout this section.

The strategy is to find a natural limsup set which covers $E(\gamma)$.

For this, we start by pointing out a property satisfied by all points in $E(\gamma)$. Heuristically, it says that every $t\in E(\gamma)$ is infinitely many times surrounded very closely by two  points which are large jumps of the Poisson point process generating $N$.

\begin{pro}
\label{pro_double}
With probability 1, the following holds: for every $t\in E(\gamma)$ and every $\e>0$ small, there exists an infinite number of integers $n\geq 0$ such that $\ds{N((t-2^{-n},t]\times [2^{-n/(\ga-\e)},1]\ge 1}$ and $\ds{N((t,t+2^{-n}]\times [2^{-n/(\ga-\e)},1]\ge 1.}$
\end{pro}

\begin{preuve}
Let $t\in E (\gamma) $. This implies that 
\begin{equation}
\label{eq1}
 \limsup_{r\to 0} \frac{\log \mu \big((\cM_t-r,\cM_t+r)\big) }{\log r} \ge \gamma\cdot\be(\cM_t) .
 \end{equation}
This equation is interpreted as the fact  that  the  time spent by the process $\cM$ in the neighborhood of $\cM_t$ cannot be too large.    The most likely way for $\mu$ to behave like this is that $\cM$  jumps into this small neighborhood of $\cM_t$ with a larger than normal jump, and quickly jumps out of that neighborhood with another big jump. This heuristic idea is made explicit by the following  computations. 

\begin{lem}
Let $0<\e< \gamma-1$.  If  $t\in E (\gamma) $, then  there exist infinitely many integers $n$ such that 
\begin{align}\label{aroundt}
|\cM_{t+2^{-n}}-\cM_t|\wedge|\cM_t-\cM_{t-2^{-n}}|\ge  2^{-\frac{n}{(\gamma-\e/4)\be(\cM_t)}}.
\end{align}
\end{lem}

\begin{preuve} \
Let us prove that $t$ satisfies 
\begin{equation}
\label{eq1911}
\limsup_{s\to 0^+}\frac{|\cM_{t+s}-\cM_t|\wedge|\cM_t-\cM_{t-s}|}{s^{1/(\gamma-\e/5)\be(\cM_t)}}\ge 1.
\end{equation}

Assume first that $\cM$ is continuous at $t$. Assume toward contradiction that for all $s>0$ sufficiently small, $|\cM_{t+s}-\cM_t| \leq s^{1/((\gamma-\e/5)\be(\cM_t))}$ or $|\cM_t-\cM_{t-s}|\leq s^{1/((\gamma-\e/5)\be(\cM_t))}$.  

If $|\cM_{t+s}-\cM_t| \leq s^{1/((\gamma-\e/5)\be(\cM_t))}$, then setting $r=\cM_{t+s}-\cM_t$,  
\begin{equation}\label{eq191}
\mu\big((\cM_t-r,\cM_t+r )\big)\geq  \mu\big((\cM_t,\cM_t+r )\big) = s \geq r^{(\gamma-\e/5)\be(\cM_t)}.  
\end{equation}
 The same holds true when $|\cM_t-\cM_{t-s}|\leq s^{1/((\gamma-\e/5)\be(\cM_t))}$. We have thus proved that  \eqref{eq191} holds for every small $r$ by continuity of $\cM$ at $t$,  this contradicts \eqref{eq1}.

\smallskip
  When $t$ is a jump time for $\cM$,  the proof goes as above using the two obvious remarks :   
$\mu(\cM_t-r,\cM_t+r) = \mu(\cM_t, \cM_t+r)$ for all small $r>0$,  and $\cM_t - \cM_{t-s} > (\cM_t-\cM_{t-})/2 $ which does not depend on $s$.

From  \eqref{eq1911}  we deduce \eqref{aroundt}.
\end{preuve} 

\smallskip

  Next technical lemma, proved in \cite{yang2015jumpdiffusion}, shows that when \eqref{aroundt} holds, there are necessarily at least two "large" jumps around (and very close to) $t$. Let us recall this lemma, adapted to our context.  
 
\begin{lem}[\cite{yang2015jumpdiffusion}]\label{lemmamarkov} Let $\wt N$ stand for the compensated Poisson measure associated with the Poisson measure $N$. There exists a constant $C$ such that  for every $\delta>1$,  for all integers $n\ge 1$ 
\begin{align*}
\pp\( \sup_{0\le s<t \le 1 \atop |s-t|\le 2^{-n}} 2^{\frac{n}{\de(\be(\cM_{t+2^{-n}})+ 2/n )}} \lba\int_s^t\int_0^{2^{-\frac{n}{\de}}} z^{\frac{1}{\be(\cM_{u-})}} \wt N(du,dz)\rba\ge 6n^2\) \le C e^{-n}.
\end{align*}
\end{lem}

\begin{remark}\label{rqmarkov}  The formula looks easier than the one in \cite{yang2015jumpdiffusion} because in our context $\cM$ is increasing.
When the function $\be$ is constant,   the term $2/n$ in the previous inequality disappears \cite{balanca2014levy}.     
\end{remark}
 
 Recall formula \eqref{defM} of $\cM$. Last Lemma allows us to control not exactly the increments of $\cM$, but the increments of the "part of $\cM$" constitued by the jumps of size less than $2^{-\frac{n}{\de}}$. It essentially entails that these ``restricted" increments over any interval of size less than $2^{-n}$ are uniformly controlled by $2^{-\frac{n}{\de(\be(\cM_{t+2^{-n}})+ 2/n )}} $ with large probability. 
  
 More precisely,  Borel-Cantelli Lemma applied to Lemma \ref{lemmamarkov} with $\delta = \gamma-\e$  yields that for    all integers $n$ greater than some $n_{\gamma-\ep}$,
\begin{align*}
\lba\int_t^{t+2^{-n}}\int_0^{2^{-\frac{n}{\gamma-\e}}} z^{\frac{1}{\be(\cM_{u-})}} \wt N(du,dz)\rba & \le 6n^2\cdot 2^{-\frac{n}{(\gamma-\e)(\be(\cM_{t+2^{-n+1}})+2/n)}}\\
& \le 2^{-\frac{n}{(\gamma-\e/2)\be(\cM_{t+2^{-n+1}}) }}.\label{secondjump}
\end{align*}
On the other hand,   for all integers $n$ greater than some other $n'_{\ga-\e}$, a direct computation gives 
\begin{align*}
\int_t^{t+2^{-n}}\int_0^{2^{-\frac{n}{\gamma-\e}}} z^{\frac{1}{\be(\cM_{u-})}} du\,\frac{dz}{z^2} & \le C2^{-n} 2^{-\frac{n}{\gamma-\e}\( \frac{1}{\be(\cM_{t+2^{-n}})}-1\) }\\
& \le 2^{-\frac{n} {(\gamma-\e)\be(\cM_{t+2^{-n}})}}.
\end{align*}
 Therefore,  for all large $n$, 
\begin{align}\nonumber
&\int_t^{t+2^{-n}}\int_0^{2^{-\frac{n}{\gamma-\e}}} z^{\frac{1}{\be(\cM_{u-})}} N(du,dz) \\  
&\le \lba\int_t^{t+2^{-n}}\int_0^{2^{-\frac{n}{\gamma-\e}}} z^{\frac{1}{\be(\cM_{u-})}} \wt N(du,dz)\rba + \int_t^{t+2^{-n}}\int_0^{2^{-\frac{n}{\gamma-\e}}} z^{\frac{1}{\be(\cM_{u-})}} du\,\frac{dz}{z^2}\nonumber \\ 
& \le 2^{-\frac{n} {(\gamma-\e/3)\be(\cM_{t+2^{-n+1}})}}.\end{align}
 Similarly, one establishes that 
\begin{align}\label{firstjump}
\int_{t-2^{-n}}^{t}\int_0^{2^{-\frac{n}{\gamma-4\e}}} z^{\frac{1}{\be(\cM_{u-})}} N(du,dz) \le 2^{-\frac{n} {(\gamma-\e/3)\be(\cM_{t+2^{-n+1}})}}.
\end{align}

Let us introduce, for every integer $n\geq 1$, the process
 $$\wt\cM^n_t= \int_0^t\int_{2^{-\frac{n}{\gamma-\e}}}^1 z^{1/\be(\cM_{u-})} N(du,dz),$$
 so that $\cM_t =\wt\cM^n_t+ \int_0^t \int_0^{2^{-\frac{n}{\gamma-\e}}}  z^{\frac{1}{\be(\cM_{u-})}} N(du,dz)$.
  
A direct estimate shows that by right-continuity of $\cM$, when $n$ becomes large, one has 
$$ 3\cdot 2^{-\frac{n} {(\gamma-\e/3)\be(\cM_{t+2^{-n+1}})}}  <    2^{-\frac{n}{(\gamma-\e/4)\be(\cM_t)}}$$

Recalling formula \eqref{defM}, the three inequalities  \eqref{aroundt},  \eqref{secondjump} and \eqref{firstjump}  imply that   for an infinite number of integers $n$
\begin{align}\label{aroundt2}
|\wt\cM^n_{t+2^{-n}}-\wt\cM^n_t|\wedge|\wt\cM^n_t-\wt\cM^n_{t-2^{-n}}|\ge  2^{-\frac{n}{(\gamma-\e/3) \be(\cM_{t+2^{-n+1}})}}\ge  2^{-\frac{n}{(\gamma-\e/2)\be(\cM_t)}} .
\end{align}
Since $\wt \cM^n$ (and $\cM$) are purely discontinuous and right continuous, this last inequality proves the existence of at least one   time $t_1^n \in (t-2^{-n},t]$ and another   time $t_2^n \in (t,t+2^{-n}]$ such that $\wt \cM^n$ (and $\cM)$ has a jump.  The desired property on the Poisson measure  $N$ follows, and Proposition \ref{pro_double} is proved. \end{preuve}

\medskip

Further, in order to find an upper bound for the dimension of $E(\gamma)$, one constructs a suitable covering of it. For $n\in\nn^*$ and $k=0,\ldots,2^n-1$, set
\begin{align*}
I_{n,k}=[k2^{-n},(k+1)2^{-n})  \ \ \mbox{ and }  \ \ \wh I_{n,k} = \bigcup_{\ell=k-1}^{k+1}I_{n,\ell}  .
\end{align*}
One introduces the collection of sets
\begin{align*}
E_n(\gamma,\e) = \llb \wh I_{n,k} : N\( \wh I_{n,k}\times\[ 2^{-\frac{n}{\gamma- \e}},1\] \) \ge 2, \quad k=0, \ldots, 2^n-1 \rrb ,
\end{align*}
 which is constituted by the intervals $\wh I_{n,k}$ containing at least two jumps for $N$ of size greater than $2^{-\frac{n}{\gamma- \e}}$.
 Finally, one  considers the limsup set
\begin{align}
\label{defegamma}
E(\gamma,\e) =  \limsup_{n\to+\infty} \ \bigcup_{\wh I\in E_n(\gamma,\e)} \wh I.
\end{align}

Proposition  \ref{pro_double} states exactly that $E(\gamma) \subset E(\gamma,\e)$.   So it is enough to find an upper bound for the Hausdorff dimension of $E(\ga,\e)$. Next lemma estimates the number of intervals contained in $E_n(\gamma,\e)$.

\begin{lem}\label{lemmadoublejump} With probability $1$, there exists a constant $C$ such that for all $n\geq 1$, 
\begin{align*}
\# E_n(\gamma,\e) \le  C n^22^{n(\frac{2}{\gamma-\e}-1)}.
\end{align*}
\end{lem} 

\begin{preuve} \ 
For a fixed "enlarged" dyadic interval $\wh I_{n,k}$,  the inclusion $\wh I_{n,k}\in E_n(\gamma,\e)$ corresponds to the event that a Poisson random variable with parameter $\mathfrak{q}_n = 3\cdot 2^{-n}\cdot 2^{\frac{n}{\gamma- \e}}$ is larger than $2$.  Since $\mathfrak{q}_n \to 0$ exponentially fast, one has   
\begin{align*}
p_n := \pp(\wh I_{n,k}\in E_n(\gamma,\e))    = C_n 2^{-n(2-\frac{2}{\gamma-\e})},
\end{align*}
where $C_n$ is a constant depending on $n$ which stays bounded away from 0 and infinity.

The events $\{\wh I_{n,3k} \in E_n(\gamma,\e)\}_{k\geq 0}$ being independent,   the random variable $\#\{k \in\{1,\ldots, \lfloor 2^n/3\rfloor: \wh I_{n,3k} \in E_n(\gamma,\e) \}$ is a Binomial random variable with parameters $(\lfloor 2^n/3\rfloor, p_n)$.
An application of Markov inequality yields
\begin{align*}
\pp\( \# \llb  k  : \wh I_{n,3k} \in E_n(\gamma,\e) \rrb \ge n^2\lfloor 2^n/3\rfloor p_n\) \le n^{-2} .
\end{align*} 
Further, Borel-Cantelli Lemma gives  that almost surely, for $n$ sufficiently large, 
$$ \# \llb  k  : \wh I_{n,3k} \in E_n(\gamma,\e) \rrb \le n^2\lfloor 2^n/3\rfloor p_n.$$
The same holds  for $\# \llb  k  : \wh I_{n,3k+1} \in E_n(\gamma,\e) \rrb $ and $\# \llb  k  : \wh I_{n,3k+2} \in E_n(\gamma,\e) \rrb $. One concludes that for all $n\in\nn$ sufficiently large, 
\begin{align*}
\# E_n(\gamma,\e) \le 3n^2\lfloor 2^n/3\rfloor p_n,
\end{align*}
which proves the claim.
\end{preuve}

\sk

Now we are in position to prove the upper bound for the Hausdorff dimension of $E(\gamma)$.
\sk

\begin{preuve}{\it \ of \eqref{reduc2} : }
Let $n_0$ be so large that the previous inequalities hold true for all integers $n\ge n_0$. 
Recalling \eqref{defegamma}, one knows that  for every $n_1\geq n_0$,  the union $ \bigcup_{n\ge n_1} \ \bigcup_{\wh I\in E_n(\gamma,\e)} \wh I $ forms a covering of $E(\gamma,\e)$, thus a covering of $E(\gamma)$.  

Let $s >  \frac{2}{\gamma-\e} -1$.   Fix $\eta>0$ and $n_1$ so large that all intervals $\wh I \in E_{n_1}(\gamma,\e)$ have a diameter less than $\eta$. Using the covering just above, one sees that the $s$-Hausdorff measure of $E(\gamma)$ is bounded above by 
\begin{align*}
\cH^{s}_\eta (E(\gamma))  & \le \sum_{n\ge n_1} \sum_{\wh I\in E_n(\gamma,\e)} |\wh I|^s   \le \sum_{n\ge n_1}   C n^22^{n(\frac{2}{\gamma-\e}-1)}  |3\cdot 2^{-n}|^s 
\end{align*}
which is a convergent series. Therefore, $\cH^{s}_\eta (E(\gamma)) = 0$ as  $n_1$ can be chosen arbitrarily large.  This leads to $\cH^{s}(E(\gamma)) = 0$ for every $s>\frac{2}{\gamma-\e}-1$. We have thus proved almost surely, 
\begin{align*}
\dimh(E(\gamma))\le \frac{2}{\gamma-\e} -1.
\end{align*}
Letting $\e\to 0$ yields the desired upper bound.
\end{preuve}

\section{Proof of Theorem \ref{theogeneral} : lower bound}
\label{seclowerbound}

The aim of this section is to get that with probability one,  \eqref{reduclow1} holds with $\cO=(0,1)$ for all  non-increasing c\`adl\`ag function $\Upsilon: [0,1]\to [1,2]$.

Recalling the notations in Theorem \ref{theogeneral}, for  simplicity, we write 
$$F(\Upsilon) = \ov E_\mu ^t\big((0,1),\Upsilon \big).$$

Let $\e>0$ and $0<b< \e$ be fixed until the end of Section \ref{sectionbeforeextend}.    We construct simultaneously  for all $\Upsilon$ with $1+2\e \le \Upsilon_{\min} \le 2-2\e $ and $\e'>0$,  a  random Cantor set  $C(\Upsilon, \e')\subset F(\Upsilon)$ with Hausdorff dimension larger than $2/({\Upsilon_{\min}}+\e')-1$. The lower bound for the Hausdorff dimension of $F(\Upsilon)$  follows.  

We explain how to extend the proof to   the functions $\Upsilon$ satisfying $\Upsilon_{\min} \in [1,2]\setminus [1+2\e,2-2\e]$  in subsection \ref{sectiongamma=1,2}.


\subsection{The time scales, and some notations}\label{subset6.1}

We aim at constructing Cantor sets inside $F(\Upsilon)$. Recalling Proposition \ref{pro_double},  some configurations for the jump points are key in this problem. More precisely, one knows that every point in $F(\Upsilon)$ is   infinitely often located in the middle of two large jumps which are really close to each other. 
So the Cantor set we are going to construct will focus on these behaviors. 

Let us define a (deterministic)  sequence of rapidly decreasing positive real numbers. First,
\begin{align*}
\begin{cases}  \eta_{1,0} = 10^{-10}, \\ 
  \eta_{1,\ell} =  \eta_{1,\ell-1}^{1+\e} \mbox{ for } 1\le \ell \le \ell_1 := \min\{ \ell\ge 1 :   \eta_{1,\ell} \le  e^{- \eta_{1,0}^{-1}} \}.  
\end{cases}
\end{align*}
By induction one defines the sequence $\{ \eta_{n,\ell} : n\in\N^*,  0\le \ell \le \ell_n \}$ as 
\begin{align*}
\begin{cases}
\eta_{n,0} = \eta_{n-1, \ell_{n-1}}, \\
 \eta_{n,\ell} = \eta_{n,\ell-1}^{1+\e}  \mbox{ for } 1\le \ell \le \ell_n := \min\{ \ell\ge 1 :  \eta_{n,\ell} \le e^{-\eta_{n,0}^{-1}}\},
\end{cases}
\end{align*}
which are our time scales.
One also sets 
\begin{eqnarray*}
 \eta_{n,\ell_n+1} & =& \eta_{n+1, 1}\\
\eta_n & = & \eta_{n,0}.
\end{eqnarray*}

The natural partition of $[0,1]$ induced by this sequence is denoted  by 
\begin{align*}\cJ_{n,\ell} =\llb J_{n,\ell,k}= [k\eta_{n,\ell}, (k+1)\eta_{n,\ell}) : k= 0,\ldots, \left\lfloor \frac{1}{\eta_{n,\ell}}\right\rfloor \rrb .
\end{align*}
 
 By convention, $J_{n,\ell, -1}=J_{n,\ell, -2} = J_{n,\ell, \[ \frac{1}{\eta_n}\] +1}=J_{n,\ell,\[ \frac{1}{\eta_n}\] +2} =\emptyset$.

One needs the enlarged intervals
 $$\wh{J}_{n,\ell, k}=\bigcup_{i=k-1}^{k+1} {J}_{n,\ell,i}.$$

\subsection{Zero jump and double jumps configuration}
Two types of jump configuration along the scales are of particular interest, since they are the key properties used to build relevant Cantor sets.   Recall that the  Poisson random measure  $N$  has  intensity $dt\otimes dz/z^2$.

\begin{definition}
\label{defcJngamma}
For any $n\in\nn^*$, $1\le \ell\le \ell_n$ and $\gamma\in[1+2\e,2-2\e]$,  define 
\begin{align}
\label{config1}
{\cJ_{n,\ell}^{z}(\gamma)} & = \llb {J}_{n,\ell, k} \in\cJ_{n,\ell} :  N\( \wh J_{n,\ell, k} \times [\eta_{n,\ell+1}^{1/\ga}, \eta_{n,\ell}^{1/\ga})  \) =0 \rrb\\
\label{config2}
\cJ^d_{n,\ell}(\ga)& =\left\{ J_{n,\ell,k}\in \cJ_{n,\ell}   : \begin{cases}
N\( {J}_{n,\ell, k-2}\times [\eta_{n,\ell}^{1/\ga}/2, \eta_{n,\ell}^{1/\ga}) \) = 1 \\
N\( {J}_{n,\ell, k+2}\times [\eta_{n,\ell}^{1/\ga}/2, \eta_{n,\ell}^{1/\ga}) \) = 1 
\end{cases}\right\}.
\end{align}
\end{definition}
\begin{remark}\label{remark0}
The superscript ``z" refers to "zero jump" while    ``d" refers to "double jump".
\end{remark}

 Let us start with   straightforward observations:
 \begin{itemize}
 \item   for $(n,\ell) \neq (n', \ell')$,  the composition (number and position of the intervals) of  $\cJ^z_{n,\ell}(\ga)$ and  $\cJ^z_{n',\ell'}(\ga)$  are independent thanks to the Poissonian nature of the measure $N$.  

 \item  The same holds true for the double jump configuration.  

 \item  Fixing $(n,\ell)$, for $|k-k'|\ge 3$,  the events $J_{n, \ell, k} \in \cJ^z_{n,\ell}(\ga)$ and  $J_{n, \ell, k'} \in \cJ^z_{n,\ell}(\ga)$  are independent.   
 
  \item The same holds for $\cJ^d_{n,\ell}(\ga)$ if one assumes that $|k-k'|\ge 5$.

 \item  For fixed $(n,\ell,k)$,  the events $J_{n, \ell, k} \in \cJ^z_{n,\ell}(\ga)$  and $J_{n, \ell, k} \in \cJ^d_{n,\ell}(\ga)$ are independent. 
\end{itemize}
 
Next probability estimate is fundamental in the sequel. 
\begin{lem}\label{chosen} For all $n\in\N^*$, $1\le \ell \le \ell_n$, $\gamma\in[1+2\e,2-2\e]$ and $J\in \cJ_{n,\ell}$, 
\begin{align}\label{eqchosen}
p_{n,\ell,\ga} &= \pp\( J \in {\cJ_{n,\ell}^z(\ga)}  \) =   \exp\( - C_{n,\ell} \eta_{n,\ell}^{(1-\frac{1+\e}{\ga})} \) \\
q_{n,\ell,\ga} &= \pp\( J \in {\cJ_{n,\ell}^d(\ga)}  \) = C_{n,\ell}'  \eta_{n,\ell}^{2-\frac{2}{\ga}}
\end{align}
where $C_{n,\ell}$, $C'_{n,\ell}$ are constants   uniformly  (with respect to $n$, $\ell$  and $\gamma$) bounded away from 0 and infinity.  
\end{lem}

\begin{preuve} \ 
The value of $p_{n,\ell, \ga}$  corresponds the probability that a Poisson random variable with parameter  $\mathfrak{p} =  3\eta_{n,\ell} \[ \eta_{n,\ell+1}^{1/\ga} - \eta_{n,\ell}^{1/\ga}\] $  equals zero, thus $p_{n,\ell,\ga} = e^{-\mathfrak{p}}$. On the other hand, each condition in \eqref{config2} relies on  the probability that a Poisson variable with parameter $\mathfrak{q} = \eta_{n,\ell}  ^{1-1/\ga}$ equals one. Hence, by independence,   $q_{n,\ell,\ga} = (e^{-\mathfrak{q}}\cdot \mathfrak{q})^2$. The result follows. 
\end{preuve}

\subsection{Random trees induced by the zero jump intervals and estimates of the number of their leaves}
In this section, one constructs for a fixed integer $n\in\N^*$ a nested collection of intervals, indexed by $0 \le \ell \le \ell_n$.  These intervals induce naturally a random tree with height $\ell_n+1$.   

One starts with any interval $J_n\in \cJ_{n}= \cJ_{n,0}$, which is the root of the tree, denoted by $\cT_{n,0} = \{J_n\}$.  Define by induction, for $1\le \ell \le \ell_n$,
\begin{align*}
\cT_{n,\ell} = \{ J\in \cJ_{n,\ell}  :   J\in \cJ^z_{n,\ell}(\ga)  \mbox{ and } J\subset \wt J \mbox{ for some }  \wt J \in \cT_{n,\ell-1}  \}.
\end{align*}
One focuses on the $J_n$-rooted random tree $\mathbb{T}_{n,\ga}(J_n) = (\cT_{n,0}, \ldots, \cT_{n,\ell_n})$.  The number of leaves of $\mathbb{T}_{n,\ga}(J_n)$, denoted by $|\mathbb{T}_{n,\ga}(J_n)|$, is the cardinality of $\cT_{n,\ell_n}$. 

\medskip

{\bf Fact:} Every point belonging to the intervals indexed by the leaves of the tree have the remarkable property that "they do not see" large jump points 
between the scales $\eta_n$ and $\eta_{n+1}$. This observation is made explicit in Lemma \ref{zerojump}.

\medskip

\begin{remark}
Observe that we dropped the index $\gamma$ in the definition of $\cT_{n,\ell}$ to ease the notations, since these sets will not re-appear in the following sections.
\end{remark}

 Our goal is to prove the following estimate on the number of leaves of $\mathbb{T}_{n,\ga}(J_n)$.

\begin{pro}\label{leaves}
With probability one, for every integer $n$ large, for  every $J_n \in\cJ_{n,0}$ and $\ga\in[1+2\e, 2-2\e]$, 
\begin{equation}
\label{eq_resume}
|\mathbb{T}_{n,\gamma}(J_n )| \ge \left\lfloor \frac{\eta_{n }}{2\eta_{n+1}} \right\rfloor.
\end{equation}
\end{pro}

The estimate of $|\mathbb{T}_{n,\gamma}(J_n)|$ is divided into several short lemmas.  

\begin{lem} For all $n\in\N^*$, $J_n\in\cJ_{n,0}$ and $\ga\in[1+2\e,2-2\e]$,  one has 
\begin{multline*}
\pp\( \#\cT_{n,1} \ge  \big(1-\log(1/\eta_{n,0})^{-2}\big)  \left\lfloor \frac{\eta_{n,0}}{\eta_{n,1}}\right\rfloor  p_{n,1,\ga} \)  \\
\ge 1 - 3 \exp\( -  \log(1/\eta_{n,0})^{-4} \left\lfloor \frac{\eta_{n,0}}{3\eta_{n,1}} \right\rfloor p_{n,1,\ga} /2\)
\end{multline*}
\end{lem}

\begin{preuve} \  
  For any $(n,\ell)$ and for $i\in\{0,1,2\}$, set
   $$\cT^i_{n,\ell} = \{ J \in \cT_{n,\ell} :  J =J_{n,\ell,3k+i} \in \cJ_{n,\ell} : k\in\N \}.$$

By independence (see the observations before Lemma \ref{chosen}),  for each $i\in\{0,1,2\}$, the number of vertices in $\cT^i_{n,1}$ is a binomial random variable with parameter $(\lfloor \eta_{n,0}/(3\eta_{n,1})\rfloor,  p_{n,1,\ga} )$.

By Chernoff inequality, for every binomial random variable $\mathfrak{X}$ with parameter $(\mathfrak{n}, \mathfrak{p})$, for any $\de\in(0,1)$,  one has 
\begin{align} \label{chernoff}
\pp(\mathfrak{X}\le (1-\de)\mathfrak{np}) \le \exp(- \de^2\mathfrak{np}/2).
\end{align}
The result follows applying \eqref{chernoff} with $\de=\log(1/\eta_{n,0})^{-2}$ for every $i$.
\end{preuve}

\begin{lem}
\label{lem444}
Set $ \ds a(n,\ell) =  (1-\log(1/\eta_{n,\ell-1})^{-2}) \left\lfloor \frac{\eta_{n,\ell-1}}{\eta_{n,\ell}} \right\rfloor p_{n,\ell,\ga}  .$

 For all $n\in\N^*$, $J_n\in\cJ_{n,0}$, $\ga\in[1-2\e, 2-2\e]$ and $2\le \ell \le \ell_n$,  a.s.
\begin{multline*}
 \pp \big( \#\cT_{n,\ell} \ge a(n,\ell)  \# \cT_{n,\ell-1}   \big) \\
\ge 1- 3\exp\( - \log(1/\eta_{n,\ell-1})^{-4} \( \#\cT_{n,\ell-1} \left\lfloor \frac{\eta_{n,\ell-1}}{3\eta_{n,\ell}} \right\rfloor \)  p_{n,\ell,\ga}/2 \) .
\end{multline*}
\end{lem}
\begin{preuve} \  Using again the remarks before Lemma \ref{chosen},  for every $i\in\{0,1,2\}$, the law of the random variable $\#\cT^i_{n,\ell}$ conditioning on $\#\cT_{n,\ell-1}$ is   binomial    with parameter $(   \#\cT_{n,\ell-1}\left\lfloor \frac{\eta_{n,\ell-1}}{3\eta_{n,\ell}} \right\rfloor , p_{n,\ell, \ga})$.  Applying \eqref{chernoff} gives the estimate. 
\end{preuve}


\begin{lem}\label{lemma6} Set
\begin{eqnarray*}
 &&b(n,\ell)= 1-  3\exp\left[ - \frac{\log(1/\eta_{n, \ell -1})^{-4}}{6(1-\log(1/\eta_{n,\ell -1})^{-2})}   \prod_{k=1}^\ell   a(n,k)  \right].
 \end{eqnarray*}
For all $n\in\N^*$, $J_n\in\cJ_{n,0}$ and $\ga\in[1-2\e, 2-2\e]$,  one has
$$
\pp\( |\mathbb{T}_{n,\gamma}(J_n)| \ge \prod_{\ell = 1}^{\ell_n}  a(n,\ell)  \)    \ge  \prod_{\ell=1}^{\ell_n} b(n,\ell) .
$$
\end{lem}

\begin{preuve} \    One has 
\begin{align*}
&\pp\( |\mathbb{T}_{n,\gamma}(J_n)| \ge \prod_{\ell = 1}^{\ell_n} a(n,\ell) \) \\
& \ge  \pp\( |\mathbb{T}_{n,\gamma}(J_n)| \ge \prod_{\ell = 1}^{\ell_n} a(n,\ell),  \ \  \#\cT_{n,\ell_n-1} \ge  \prod_{\ell = 1}^{\ell_n-1} a(n,\ell) \) \\
& \ge \pp\( |\mathbb{T}_{n,\gamma}(J_n)| \ge a(n,\ell_n)  \#\cT_{n,\ell_n-1} ,  \ \  \#\cT_{n,\ell_n-1} \ge   \prod_{\ell = 1}^{\ell_n-1} a(n,\ell)  \) .
 \end{align*}
Conditioning on $\#\cT_{n,\ell_n-1}$,  and using Lemma \ref{lem444} with $\ell=\ell_n$, this probability is greater than
\begin{align*}&  \E\[ \E\[  \left. 1- 3\exp\( - \log(1/\eta_{n,\ell_n-1})^{-4} \#\cT_{n,\ell_n-1} \left\lfloor \frac{\eta_{n,\ell_n-1}}{3\eta_{n,\ell_n}} \right\rfloor  p_{n,\ell_n,\ga}/2 \)    \right\vert \#\cT_{n,\ell_n-1}  \] \right.\\
 &\hspace{80mm}   \times  \indiq_{\#\cT_{n,\ell_n-1} \ge  \prod_{\ell = 1}^{\ell_n-1} a(n,\ell) }  \Big] \\
& \ge\E\[ \E\[  \left. 1- 3\exp\( -\frac{ \log(1/\eta_{n,\ell_n-1})^{-4} }{ 6(1-\log(1/\eta_{n,\ell_n -1})^{-2})} a(n,\ell_n) \#\cT_{n,\ell_n-1}  \)    \right\vert \#\cT_{n,\ell_n-1}  \] \right.\\
 &\hspace{80mm}  \times \indiq_{\#\cT_{n,\ell_n-1} \ge  \prod_{\ell = 1}^{\ell_n-1} a(n,\ell) }  \Big] \\
&\ge b(n,\ell_n)  \ \pp\( \#\cT_{n,\ell_n-1} \ge  \prod_{\ell = 1}^{\ell_n-1} a(n,\ell) \)
\end{align*}
Iterating this computation yields the desired inequality. 
\end{preuve}


\medskip

We are now in position to prove Proposition \ref{leaves}.

\medskip

\begin{preuve} \  
We are going to prove  the following lemma:
\begin{lem}
\label{lemestimate}
For some constant $c_1$,   for all $n\in\N^*$ large enough, for every $J_n\in\cJ_{n,0}$ and $\ga\in[1+2\e, 2-2\e]$,  one has
\begin{align}
\label{eq_resume2}
\pp\(  |\mathbb{T}_{n,\gamma}(J_n)| \ge \left\lfloor \frac{\eta_{n }}{2\eta_{n+1}} \right\rfloor \) 
\ge \exp\( -c_1 \exp\(  -\eta_{n}^{-\e/2}   \)  \) .
\end{align}
\end{lem}

\begin{preuve} \  
One  combines the estimates in Lemma \ref{chosen} and  Lemma \ref{lemma6}.  
Let us first estimate $\prod_{\ell=1}^{\ell'} a(n,\ell)$ for $2\leq \ell'\leq \ell_n$. Observe  that  for large $n$, one has
\begin{align*}
\ell_n &= \log\( \frac{\log(1/\eta_{n,\ell_n})}{\log(1/\eta_{n,0})} \) / \log(1+\e) \le \log\log(1/\eta_{n,\ell_n})/\log(1+\e) \\ 
&\le \log\( \frac{1+\e}{\eta_{n,0}} \) /\log(1+\e) \le 2 \log(1/\eta_{n,0}).
\end{align*}
Thus for $n$ large enough,
\begin{align}
&\prod_{\ell=1}^{\ell' } (1-\log(1/\eta_{n,\ell-1})^{-2}) \nonumber \\
&= \exp  \( \sum_{\ell=1}^{\ell' } \log\Big( 1-\log(1/\eta_{n,\ell-1})^{-2}\Big) \) \ge \exp\( -2 \sum_{\ell=1}^{\ell' } \log(1/\eta_{n,\ell-1})^{-2} \)  \nonumber \\
&\ge  \exp\( -2\ell'  \log(1/\eta_{n,0})^{-2} \) \ge \exp\( -4\log(1/\eta_{n,0})^{-1} \)  \ge 1/\sqrt{2}. \label{eq1lem7}
\end{align}
Using the rapid decay of $(\eta_{n,\ell})$  to zero and the uniform boundedness of $C_{n,\ell}$, one can find a    constant $c_0>0$ such that for all $n$ large enough,
\begin{equation}\label{eq2lem7}
\prod_{\ell=1}^{\ell'} p_{n,\ell,\ga} = \exp\( -\sum_{\ell = 1}^{\ell'} C_{n,\ell}  \eta_{n,\ell}^{1-\frac{1+\e}{\ga}} \) \ge \exp\(  -c_0 \eta_{n,1}^{1-\frac{1+\e}{\ga}} \) \ge 1/\sqrt{2}. 
\end{equation}
Combing \eqref{eq1lem7} and \eqref{eq2lem7},  one concludes that for all large $n$   
\begin{equation}
\label{minutile}
\prod_{\ell=1}^{\ell'}  a(n,\ell)\ge \left\lfloor \frac{\eta_{n,0}}{2\eta_{n,\ell'}}\right\rfloor, \mbox{    and in particular, } \ds \prod_{\ell=1}^{\ell_n}  a(n,\ell)\ge   \left\lfloor \frac{\eta_{n }}{2\eta_{n+1}}\right\rfloor  .
\end{equation}

Now we estimate  the other product $\prod_{\ell=1}^{\ell_n} b(n,\ell) $. Using \eqref{eq1lem7}, \eqref{eq2lem7} and \eqref{minutile},   there exists  a     constant $c_1>0$ such that for every large $n$
\begin{align*}
\prod_{\ell=1}^{\ell_n} b(n,\ell) &\ge \prod_{\ell=1}^{\ell_n}  \llb1-3 \exp\( -\frac{\log(1/\eta_{n,\ell-1})^{-4}}{12} \left\lfloor \frac{\eta_{n,0}}{\eta_{n,\ell}}\right\rfloor \) \rrb \\
& =  \exp  \llb   \sum_{\ell=1}^{\ell_n}  \log \( 1-3 \exp\( -\frac{\log(1/\eta_{n,\ell-1})^{-4}}{12} \left\lfloor \frac{\eta_{n,0}}{\eta_{n,\ell}}\right\rfloor \)  \) \rrb \\ 
&\ge \exp  \llb  -6  \sum_{\ell=1}^{\ell_n}   \exp\( -\frac{\log(1/\eta_{n,\ell-1})^{-4}}{12} \left\lfloor \frac{\eta_{n,0}}{\eta_{n,\ell}}\right\rfloor \)  \rrb \\ 
&\ge  \exp \llb -c_1 \exp\(  -\frac{\log(1/\eta_{n,0})^{-4}}{12} \left\lfloor \frac{\eta_{n,0}}{\eta_{n,1}}\right\rfloor \) \rrb \\
&\ge \exp\( -c_1 \exp\(  -\eta_{n,0}^{-\e/2} \) \)  = \exp\( -c_1 \exp\(  -\eta_{n}^{-\e/2} \) \) .
\end{align*}
where   the fast decay rate of $(\eta_{n,\ell})$   to zero has been used for the third inequality.

These last equations prove exactly \eqref{eq_resume2}.
\end{preuve}

Finally, to prove Proposition \ref{leaves}, since the cardinality of $\mathcal{J}_n $ is less than $\eta_n^{-1}$, 
\begin{eqnarray*}
\mathbb{P} \left(\exists J_n\in \mathcal{J}_n:  |\mathbb{T}_{n,\gamma}(J_n)| < \left\lfloor \frac{\eta_{n }}{2\eta_{n+1}} \right\rfloor \right) \leq \eta_n^{-1} \Big( 1-\exp \( -c_1 \exp \(  -\eta_{n}^{-\e/2}   \)  \) \Big).
\end{eqnarray*}
Using the fast decay of $\eta_n$ to zero, this is the general term of a convergent series, and the Borel-Cantelli lemma gives the result.
\end{preuve}

\begin{remark} Essentially, one needs to keep in mind that  the number of leaves of the random tree $\mathbb{T}_{n,\gamma}(J_n)$ is the total number of intervals of $\cJ_{n+1}$ inside $J_n$,  up to a constant factor 1/2.  
\end{remark}

One finishes this section by proving that every point belonging to  a leaf of $\mathbb{T}_{n,\gamma}(J)$ "is not close"  to large jumps.
\begin{lem}
\label{zerojump}
Let $J\in\cJ_n$ and  $r\in [\eta_{n+1}, \eta_n)$.   Assume that $\mathbb{T}_{n,\ga}(J)$ is not empty.  Then for each $t\in \mathbb{T}_{n,\ga}(J)$,   $$N(B(t,r)\times [r^{1/\Upsilon^n_{J_{n,0}(t)}}, \eta_{n}^{1/\Upsilon^n_{J_{n,0}(t)}}]  ) =0.$$
\end{lem}
\begin{preuve}\ 
For each $t\in \mathbb{T}_{n,\ga}(J)$,  denote by $J_{n,\ell}(t)$ the unique interval such that $t\in J_{n, \ell}(t)$ for all $0\le \ell\le \ell_n$.  Denote by $\ell_0$ the unique integer such that $\eta_{n, \ell_0+1}\le r < \eta_{n,\ell_0}$.     By construction of the random tree $\mathbb{T}_{n, \ga}(J)$,  one has 
$$N(B(t,r)\times [r^{1/\Upsilon^n_{J_{n,0}(t)}}, \eta_{n,\ell_0}^{1/\Upsilon^n_{J_{n,0}(t)}}]  ) \le N(\wh J_{n,\ell_0}(t)\times  [\eta_{n,\ell_0+1}^{1/\Upsilon^n_{J_{n,0}(t)}}, \eta_{n,\ell_0}^{1/\Upsilon^n_{J_{n,0}(t)}}])=0.$$
Further,  all ancestor interval of $J_{n,\ell_0}(t)$ has no large jumps around,  in particular, 
$$N(\wh J_{n,\ell_0}(t)\times  [\eta_{n,\ell_0}^{1/\Upsilon^n_{J_{n,0}(t)}}, \eta_{n}^{1/\Upsilon^n_{J_{n,0}(t)}}])=0$$
Combining these estimates yields the result. 
\end{preuve}

\subsection{Double jumps configuration around the leaves, and key lemma}

In the previous section, we have seen that the "zero jump" configuration is quite frequent. The aim here is to estimate the number of  intervals with "double jumps" amongst the leaves of the trees. To this end,  we introduce further some notations. Set
 $$M_{n}(\ga) = \eta_{n+1}^{1-2/(\ga+ 3\cdot 2^{-n-1})}\eta_{n}^3.$$ 

\sk
\begin{definition}
Let $J_0\in\cJ_{n,0}$ and $\mathbb{T}_{n,\ga}(J_0)$ be  the random tree defined in last subsection. Consider its leaves that we denote  $\{J'_i\}_{i=1,\ldots, |\mathbb{T}_{n,\ga}(J_0)|}$, which are intervals of length $\eta_{n+1}$.

The   families  $\{\cF(J_0,\gamma, m)\}_{m=1,..., \lfloor M_{n}(\ga)/2\rfloor }$  are defined as the following disjoint subfamilies of  $\{J'_i\}_{i=1,\ldots, |\mathbb{T}_{n,\ga}(J_0)|}$ : 
$$\cF(J_0,\ga,m) = \left\{ J'_{5k + 2m \left\lfloor \frac{|\mathbb{T}_{n,\ga}(J_0)|}{M_{n}(\ga) }\right\rfloor }: k\in \left\{0,..., \left \lfloor  \frac{|\mathbb{T}_{n,\ga}(J_0)|}{5M_{n}(\ga) } \right \rfloor \right\}\right\}.$$
\end{definition}

Hence, two families $\cF(J_0,\ga,m) $ and $\cF(J_0,\ga, m')$ are disjoint and separated by a distance equivalent to $  \lfloor |\mathbb{T}_{n,\ga}(J_0)|/{M_{n}(\ga)}\rfloor\eta_{n+1}$, and the intervals belonging to the same $\cF(J_0,\ga, m) $ are separated by the distance at least $4\eta_{n+1}$.

Finally, denote by 
 $$\cD_n = \{ k2^{-n} : k\in\mathbb{Z}, n\in\nn^*\}$$
 the  $n$-th generation dyadic numbers.  One is ready to prove the key lemma.

\begin{lem}\label{numberbiggerthan} The following holds with probability $1$: there exists a (random) integer $n_0$ such that for all $n\geq n_0$,  for every  ${J}\in\cJ_{n}$, every $\gamma\in\cD_{n}\cap[1+2\e,2-2\e]$, every $\mathfrak{a} \in \{0,1,2,3\}$,  each family   $\{\cF(J,\ga, m)\}_{m=1,..., \lfloor M_{n}(\ga)/2\rfloor }$  contains  at least one interval belonging to $ \cJ^d_{n+1}(\gamma + \mathfrak{a}\cdot 2^{-(n+1)})$.
\end{lem}

\

Remark that the intervals belonging to $\cF(J,\ga, m)$ come also from the tree $\mathbb{T}_{n,\gamma}(J)$ associated with $J$, so they also enjoy the "zero jump" property.

 \
 
\begin{preuve} \ 
Fix  $n$ a positive integer, $J\in\cJ_{n }$,  $\ga\in \cD_{n}\cap[1+2\e, 2-2\e]$, $\mathfrak{a}\in \{0,1,2, 3\}$.    

Recall that $\cJ_{n }=\cJ_{n ,0}$ and $\cJ_{n +1}=\cJ_{n ,\ell_n}$ with the notations of the previous sections.

By Lemma  \ref{chosen} and the observations made before this Lemma, there exists a positive finite constant $c_2$ such that for all $n$ large
\begin{align*}
&\pp\(  \exists m, \ \forall J'\in \cF(J, \ga + \mathfrak{a}  2^{-n-1}, m),  \ J'\not\in \cJ^d_{n, \ell_n}(\ga + \mathfrak{a}  2^{-n-1}) \right. \\
& \hspace{50mm} \left. \ \left\vert  \ |\mathbb{T}_{n,\ga+\mathfrak{a}2^{-n-1}}(J)| \) \ge \left\lfloor \frac{\eta_{n}}{2\eta_{n+1}} \right\rfloor    \) \\
&\le \left\lfloor \frac{M_n(\ga)}{2} \right\rfloor  \( 1- q_{n,\ell_n, \ga + \mathfrak{a}\cdot 2^{-n-1}} \) ^{\frac{1}{M_n(\ga)} \cdot \left\lfloor\frac{\eta_{n}}{10\eta_{n+1}} \right\rfloor } \\
&\le  \eta_{n+1}^{1-\frac{2}{\ga + 3\cdot 2^{-n-1}} }\eta_{n}^3   \( 1- c_2 \eta_{n+1}^{2-\frac{2}{\ga + \mathfrak{a}  2^{-n-1}}}\) ^{ \frac{\eta_{n}}{10\eta_{n+1}^{2- \frac{2}{\ga + 3\cdot 2^{-n-1}}} \eta_{n}^3  }  }  \\
&\le  \eta_{n+1}^{1-\frac{2}{\ga + 3\cdot 2^{-n-1} }}\eta_{n}^3  \exp\( - c_2\eta_{n}^{-2}/10  \)
\end{align*}
Remark that $\eta_{n+1} \le e^{-\eta_{n}^{-1}} \le \eta_{n,\ell_n-1}$ implies $\log(1/\eta_{n+1})\le (1+\e)\eta_{n}^{-1}$.  The above probability  is thus  bounded by above by
\begin{align*}
\eta_{n}^3    \exp\(  \( \frac{2}{\ga + 3\cdot 2^{-n-1}}-1\) (1+\e)\eta_{n}^{-1} - \frac{c_2}{10} \eta_{n}^{-2} \) \le \eta_{n}^3.
\end{align*}

On the other hand,  by Lemma \ref{lemestimate} one has 
\begin{multline*}
\pp\( |\mathbb{T}_{n,\ga + \mathfrak{a}  2^{-n-1}}(J)| \le  \left\lfloor \frac{\eta_{n}}{2\eta_{n+1} }  \right\rfloor \)  \\ 
\le 1-  \exp\( -c_1 \exp(-\eta_{n}^{-\e/2}) \) \le  2c_1 \exp(-\eta_{n}^{-\e/2}).
\end{multline*}
Thus, 
$\ds \pp\(  \exists m, \ \forall J'\in \cF(J, \ga, m),  \ J'\not\in \cJ^d_{n+1}(\ga + \mathfrak{a}  2^{-n-1})  \)  \le   2\eta_{n}^{3} $.  One deduces that  
\begin{align*}
   &\pp\big( \exists J \in \cJ_{n }, \, \exists\,m, \, \forall J' \in \cF(J,\ga,m), \, J'\notin  \cJ^d_{n+1}(\gamma+\mathfrak{a}  2^{-n-1})\big )\\
   &\leq  \eta_{n}^{-1} \cdot 2\eta_{n}^{3}   = 2 \eta_{n}^{2}.
\end{align*}
There are less than $2^n$ possible choices for $\gamma$, and 4 choices for $\mathfrak{a}$. Hence,
\begin{align*}
   &\pp\big( \exists\, \gamma, \exists\, \mathfrak{a},  \exists \, J \in \cJ_{n }, \, \exists\,m,  \forall J' \in \cF(J,\ga, m), \ J'\notin  \cJ^d_{n+1}(\gamma+\mathfrak{a}  2^{-n-1})\big ) \\
   &\leq 2^{n+3}\eta_{n}^{2},
\end{align*}
which is the general term of a convergent series.  An application of  Borel-Cantelli Lemma entails the result.  
\end{preuve}

\subsection{Construction of the Cantor sets}

We are ready to construct the families of Cantor sets  $\{C(\Upsilon,\e')$ associated with c\`adl\`ag non-increasing  functions $\Upsilon: \zu \to [1+2\e, 2-2\e]$, where $\e' $ is any positive rational parameter. These sets are constituted by points which only see those double jump configurations studied before, and their Hausdorff  dimension satisfies $\dimh C(\Upsilon,\e') \geq \frac{2}{\Upsilon_{\min}+2\ep'} -1.$

\mk

{\bf Step 1 (localization).} For each $\Upsilon$ as above and $\e'>0$, there exist $t_{\e'}\in(0,1)$, $r_{\e'}>0$ such that $\forall\, t\in [t_{\e'}-r_{\e'}, t_{\e'}+r_{\e'}]$,  we have $\Upsilon(t) < \Upsilon_{\min} +\e'$. 

Let $n_0$ be the random integer  obtained in Lemma \ref{numberbiggerthan}. We  assume that   $n_0$   is so  large  that the conclusions \eqref {eq_resume} of Proposition \ref{leaves}  hold, and also that 
$$2\e'/\eta_{n_0} > K_\Upsilon\cdot 2^{n_0}, \ \ \mbox{ where } \ K_\Upsilon = |\Upsilon(1-)-\Upsilon(0)|<+\infty.$$

For every interval $J$, let $\mbox{Osc}_\Upsilon({J}) = \sup_{t\in{J}}\Upsilon(t) -\inf_{t\in{J}}\Upsilon(t)$ be the oscillation of $\Upsilon$ over $J$. 
By the monotonicity of $\Upsilon$,  for each $n\ge n_0$ one has
\begin{equation}\label{restriction} 
\#\{{J} \in\cJ_{n} : {J} \subset[t_{\e'}-r_{\e'}, t_{\e'}+r_{\e'}]  \mbox{ and } \mbox{Osc}_\Upsilon({J})\ge 2^{-n} \} \le K_\Upsilon\cdot2^{n}
\end{equation}

\mk

{\bf Step 2 (Initialization of  the Cantor set). } One chooses arbitrarily one interval ${{J}^{{n_0}}} \in\cJ_{n_0 }$ contained in $ [t_{\e'}-r_{\e'}, t_{\e'}+r_{\e'}]$ such that 
\begin{equation}\label{restriction0}
 \mbox{Osc}_\Upsilon({{J}^{{n_0}}})<2^{-n_0}.
\end{equation}
Set   the    generation "zero" of the Cantor set  as $\cC_{{n_0}}(\Upsilon,\e') ={{J}^{{n_0}}}$.

Simultaneously, we build a measure $\nu_{n_0}$ by setting $\nu_{n_0}({{J}^{{n_0}}}) =1$, and $\nu_{{n_0}}$ is uniformly distributed on ${{J}^{{n_0}}}$.

\mk

{\bf Step 3 (Next generation of the Cantor set). }
 
 Let us introduce the following notation: for each $n\in\nn^*$ and ${J}\in\cJ_{n }$, set 
 $$\Upsilon_{{J}}^n = \max\left(\cD_n\cap[1,\inf_{t\in{J}}\Upsilon(t)]\right).$$
 
   We explain how to get the second generation of intervals  $\cC_{{n_0+1}}(\Upsilon,\e')$ of the Cantor set.

The oscillation restriction \eqref{restriction0} for $\Upsilon$ on ${J}^{{n_0}}$ implies that  for every ${J}\in\cJ_{n_0+1 }$ contained in ${J}^{{n_0}}$, the quantity $\Upsilon_{{J}}^{n_0+1}$ takes necessarily one of the    four  values  $\{ \Upsilon_{{J}^{{n_0}}}^{n_0} + \mathfrak{a}  2^{-n_0-1} : \mathfrak{a}= 0,1,2, 3\}$.

Moreover, applying Lemma \ref{numberbiggerthan} to ${J}^{{n_0}}$, one obtains that   for   each  $\mathfrak{a}\in\{ 0,1,2, 3\}$, each subfamily $\{\cF({J}^{{n_0}}, \Upsilon^{n_0}_{{J}^{{n_0}}} , m)\}_{m=1,..., \lfloor M_{n_0}( \Upsilon_{{J}^{{n_0}}}^{n_0} )/2\rfloor }$  contains  at least one interval $J$ belonging to $  \cJ^d_{n_0+1 }( \Upsilon_{{J}^{{n_0}}}^{n_0}  + \mathfrak{a}  2^{-n_0-1})$.

Recalling that $\Upsilon $ is non-increasing, the quantities $\Upsilon_{{J}}^{n_0+1}$ are also non-increasing when $J$ ranges from left to right.  Since there are $ \lfloor M_{n_0}( \Upsilon_{{J}^{{n_0}}}^{n_0} )/2\rfloor$ disjoint families $\{\cF(J_{0}, \Upsilon^{n_0}_{{J}^{{n_0}}} , m)\}$  which are organized  in increasing order,  one deduces that there is $\mathfrak{a} \in\{0,1,2,3\}$ such that   there exist $\lfloor M_{n_0} ( \Upsilon_{{J}^{{n_0}}}^{n_0} )/8\rfloor$   different integers  $m \in \{1,...,  \lfloor M_{n_0}( \Upsilon_{{J}^{{n_0}}}^{n_0} )/2\rfloor\}$ for which the family $\cF(J_{0},\Upsilon^{n_0}_{J^{n_0}}, m)$ contains (at least) one interval ${J}$ satisfying $\Upsilon_{{J}}^{n_0+1}=  \Upsilon_{{J}^{{n_0}}}^{n_0}  + \mathfrak{a}  2^{-n_0-1} $  and $J \in  \cJ^d_{n_0+1} (\Upsilon_{{J}}^{n_0+1})$.

Remark that 
$$2^{n_0+1}\ll  \lfloor M_{n_0}( \Upsilon_{{J}^{{n_0}}}^{n_0} )/8\rfloor/2$$
 where we used that $\Upsilon_{\min}<2-2\e$ and ${J}^{{n_0}} \subset [t_{\e'}-r_{\e'}, t_{\e'}+r_{\e'}]$. Then, applying \eqref{restriction} for $n=n_0+1$,  one can choose  the first $\lfloor M_{n_0}(\Upsilon^{n_0}_{{J}^{{n_0}} }  )/16 \rfloor  $ intervals $J$ which satisfy $\mbox{Osc}_\Upsilon(J) < 2^{-n_0-1}$ among those already  selected in the last paragraph.  
 
 \sk
 
Finally,  $C_{{n_0}}(\Upsilon,\e')$ is the union of these intervals, which are called the basic intervals of generation ${n_0}+1$. Observe that these intervals are separated by a distance larger than $ \eta_{n_0}/ (2{M_{n_0}(\Upsilon^{n_0}_{{J}^{{n_0}} })} )$ (thanks to Borel-Cantelli applied to Lemma \ref{leaves}),  and they all have their length equal to $  \eta_{n_0+1 }$.

\mk

Simultaneously, one defines a refinement $\nu_{{n_0}+1}$ of the measure $\nu_{0}$  by setting for every $J^{n_0+1}$ basic interval of $C_{{n_0}+1}(\Upsilon,\e')$
$$ \nu_{{n_0}+1}(J^{n_0+1}) = \nu _{{n_0}}({J}^{{n_0}} ) \frac {1}{\lfloor M_{n_0}(\Upsilon^{n_0}_{{J}^{{n_0}} }  )/16 \rfloor },$$
and by saying that $ \nu_{{n_0}+1}$ is uniformly distributed inside each $J^{n_0+1}$.

\mk

{\bf Step 4: Induction of the construction of the Cantor set:}
   
Assume that for every $i=n_0,n_0+1,..., n_0+n$, the generation $\cC_{i}(\Upsilon,\e')  $ has been constructed and satisfies the following:
\begin{enumerate}
\sk\item
$\cC_{i}(\Upsilon,\e')  $ is constituted by a finite number of basic disjoint  intervals $J^{ i}$ belonging to $\cJ_{i }$, 

\sk\item
for every $i= n_0+1,..., n_0+n$, each basic interval $J^{i}\in \cC_{i}(\Upsilon,\e') $ is included in a unique basic interval $J^{i-1} \in  \cC_{ i-1}(\Upsilon,\e') $.

\sk\item
for every $i=n_0,n_0+1,..., n_0+n-1$, each basic interval $J^{i}\in \cC_{ i}(\Upsilon,\e') $  contains  $\lfloor M_{  i}(\Upsilon^{  i}_{{J}^{i}}  )/16 \rfloor  $ intervals $J^{i+1} \in \cC_{ i+1}(\Upsilon,\e') $. These intervals are separated by a distance at least equal to $\eta_{ i}/ (2  {M_{ i}(\Upsilon^{ i}_{{J^i}})})$.  Moreover, each $J^{i+1}$ belongs to $ \cJ^d_{i+1 }(\Upsilon_{J^{i+1}}^{ i+1})$.

\sk\item
Each basic interval $J^i$ of $\cC_{i}(\Upsilon,\e')  $  satisfies   $\mbox{Osc}_\Upsilon(J^i) \le  2^{-i}$.  
 
 \sk\item
 for every $i= n_0+1,..., n_0+n$, $\nu_{ i}$ is a measure supported by the basic intervals $J^i$ of $\cC_{ i}(\Upsilon,\e')  $, and if $J^{i-1}$ is the unique interval in $\cC_{ i-1}(\Upsilon,\e') $ such that $J^i \subset J^{i-1}$, then
 \begin{equation}
 \label{defnu}
  \nu_{ i}(J^i) = \nu _{ i-1}({J}^{i-1} ) \frac {1}{\lfloor M_{ i-1}(\Upsilon^{ i-1}_{{J}^{i-1}}  )/16 \rfloor }
  \end{equation}
and $\nu_{ i}$ is uniformly distributed inside each $J^i$.
\end{enumerate}

 We are now able to complete the induction.
 
For any basic interval $J^{ n} \in \cC_{ n}(\Upsilon,\e')$, applying the same method as in step 3, one finds $\lfloor M_{n}(\Upsilon^{ n}_{J^{ n}  }  )/16 \rfloor  $ intervals $J^{ n+1} \in \cJ^d_{ n+1  }(\Upsilon_{{J}^{ n+1}}^{ n+1})$,   also satisfying  $\mbox{Osc}_\Upsilon(J^{ n+1} ) < 2^{ -(n+1)}$.
 
Then,  $C_{ n+1}(\Upsilon,\e')$ is the union of these intervals, which constitue  the basic intervals of generation $  n+1 $. By construction, these basic intervals $J^{ n+1}$  are separated by at least $\eta_{ n}/ (2 M_{ n}(\Upsilon^{ n}_{J^{ n}}))$, (where $J^{ n}$ is the ``parent" interval of $J^{ n+1}$, i.e. the unique  basic interval  in  $ \cC_{ n}(\Upsilon,\e')$ such that $J^{ n+1}\subset J^{ n}$), and they all have their length equal to $  \eta_{ n+1}$.

\mk

Simultaneously, the refinement $\nu_{ n+1}$ of the measure $\nu_{ n}$  is defined by setting, for every $J^{ n}$ is the "parent" interval of $J^{ n+1}$,
$$ \nu_{ n+1}(J^{ n+1}) = \nu _{ n}(J^{ n}) \frac {1}{\lfloor M_{ n}(\Upsilon^{  n }_{{J}^{ n}}  )/16 \rfloor } ,$$
and by saying that $ \nu_{ n+1}$ is uniformly distributed inside each $J^{ n+1}$.

\begin{pro}
The Cantor set $C(\Upsilon,\e')$ is defined as
$$\cC(\Upsilon,\e') = \bigcap_{n\ge n_0} \bigcup_{{J}\in \cC_n(\Upsilon,\e')} {J}.$$
There exists a unique Borel probability measure $\nu _{\Upsilon,\e'}$ supported exactly by $\cC(\Upsilon,\e')$ such that  for all $n\ge n_0$,  the measure $\nu _{\Upsilon,\e'}$ restricted to the $\sigma$-algebra generated by $\{  \cJ_{ k}:   n_0\le k \le n  \}$ is   $\nu_{ n}$.
\end{pro}
The proof is immediate, since the  step 4. of the construction ensures that the measure is a   well-defined additive set function with total mass 1 on the algebra $\{  \cJ_{ n}:     n \geq n_0  \}$ which generates the Borel $\sigma$-algebra, thus extends to a unique probability measure on Borel sets.

\sk
Observe that the construction of the family of Cantor sets depends only on Lemma \ref{numberbiggerthan}, which holds with probability one simultaneously for all functions $\Upsilon$, as desired.

\mk

\subsection{Properties of the Cantor sets} 

The following proposition is key, since it shows that our construction guarantees that we have built points in $F(\Upsilon)$.

\begin{pro}
Almost surely, for every non-increasing c\`adl\`ag function $\Upsilon: [0,1]\to [1+2\ep,2-2\ep]$ and for every small $\ep'>0$, 
$$\cC(\Upsilon,\e') \setminus \big(S(\cM)\cup S(\Upsilon)\big)\subset F(\Upsilon).$$
\end{pro}

\begin{preuve} \ 
Suppose that $t\in[0,1]\cap \cC(\Upsilon,\e')$ is a continuous time for $\cM$ and $\Upsilon$. One wants to prove that $\ov \dim(\mu,\cM_t) = \Upsilon(t)\be(\cM_t)  $.
 
\sk 

We start by bounding   $\ov \dim(\mu,\cM_t) $ from below.

%

By construction, for every $n\geq n_0$,  $t$ is covered by an interval ${J}^n \in \cJ^d_{n}(\Upsilon_{{J}^n}^{n})$, a basic interval in $ \cC_{n }(\Upsilon,\e')$. Since $J^n  \in \cJ^d_{n}(\Upsilon_{{J}^{n}}^{n}) $, property
  \eqref{config2}   entails that  $t \in J^n$ is surrounded by two  jumps of the Poisson point process   located at $t^1_n$ and $t^2_n$ whose size belong to $\in [\eta_{n}^{1/ \Upsilon^n_{J^n}}/2, \eta_{n}^{1/ \Upsilon^n_{J^n}}] $,  and   whose mutual distance is at least $3\eta_{n}$, and at most $5 \eta_{n}$. 
  
  The process $\cM$ jumps at $t^1_n$ and $t^2_n$, with jump  size $\eta_{n}^{1/ (\Upsilon^n_{J^n}\be(\cM_{t^1_n-})) }$,   and  $\eta_{n}^{1/ (\Upsilon^n_{J^n}\be(\cM_{t^2_n-})) }$.  Since the process $\cM$ is increasing, if $J^n $ is written $ [k\eta_n, (k+1)\eta_n)$, both size of jumps at $t^1_n$ and $t^2_n$ are bounded by  below by 
  $$  r_n= \eta_{n}^{\frac{1}{\Upsilon_{{J^n}}^n \be(\cM_{(k-3)\eta_{n } })} }. $$
    Hence, $  \mu\left( B\left(\cM_t , r_n\right)\right)  \leq  5\eta_{n}$. Applying this when  $n$ becomes large, one gets
\begin{align*}
\ov \dim(\mu,\cM_t)&\ge \limsup_{n\to+\infty} \frac{\log \mu\left( B\left( \cM_t ,r_n\right)\right)   }{\log r_n}  \geq  \limsup_{n\to +\infty} \frac{  \log 5\eta_n }{\log r_n } =  \Upsilon(t)\be(\cM_t) ,
\end{align*}
where we used the continuity of $\cM$ and $\Upsilon$ at $t$.    
 
 \mk

The rest of the proof is dedicated to  prove the converse inequality, i.e.  $\overline{\dim}(\mu,\cM_t)\le \Upsilon(t)\cdot\be(\cM_t)$, which is more delicate.

\medskip

Let $\e_1>0$ be small. Thanks to the continuity of $\cM$ at $t$,  there exists $r_0 >0$ such that $r_0\le \eta_{n_0 }$ and
\begin{equation*}
N([t-r_0,t+r_0]\times [\eta_{n_0,}^{1/\Upsilon^{n_0}_{J^{n_0}}}, 1]) = 0, 
\end{equation*}
where $J^{n_0}$ is the unique interval of $C_{n_0}(\Upsilon, \e')$ that contains $t$. 

 Now for any $0<r<r_0/3$, there exists a unique integer $n \geq n_0$   such that $\eta_{n +1} \le r <\eta_{n }$. Let us call $J_{n }(t)$  and $J_{n+1}(t)$ the unique intervals of $\cJ_{n }$ and  $\cJ_{n+1}$ that contain $t$. 
 
 By construction of the random tree $\mathbb{T}_n(J_{n }(t))$,  there is no large jump around $t$. More precisely, by Lemma \ref{zerojump},    
 $$N\( B(t,r)\times [r^{1/\Upsilon^n_{J_{n}(t)}}, \eta_{n}^{1/\Upsilon^n_{J_{n}(t)}}]  \)  =0.$$
Applying same argument as in Lemma  \ref{zerojump} to scales between $n_0$ and $n$, together with the fact that the sequence $n\mapsto \Upsilon^{n+1}_{J_{n+1}(t)}$ is increasing, yields that 
$$N\( B(t,r) \times [\eta_{n }^{1/\Upsilon^{n}_{J_{n}(t)}}, \eta_{n_0 }^{1/\Upsilon^{n_0}_{J^{n_0}}}]\) = 0.$$
One deduces that the increment of $\cM$ between $t-r$ and $t+r$ has the form
\begin{align*}
\cM_{t+r} - \cM_{t-r} = \int_{t-r}^{t+r} \int_{0}^{r^{1/\Upsilon^{n}_{J_{n }(t)}}} z^{1/\be(\cM_{s-})} N(ds, dz).
\end{align*} 
Denote by $m$ the unique integer such that $2^{-m-1}\le 2r< 2^{-m}$. One has 
\begin{align*}
\cM_{t+r} - \cM_{t-r}  &\le \lba  \int_{t-r}^{t+r} \int_0^{2^{-m/\Upsilon^{n}_{J_{n }(t)} }}    z^{1/\be(\cM_{s-})} \wt N(ds,dz) \rba \\
  &+ \int_{t-r}^{t+r}  \int_0^{2^{-m/\Upsilon^{n}_{J_{n }(t)} }}    z^{1/\be(\cM_{s-})} \,\frac{dz}{z^2}\, ds := A_1(r)+A_2(r).
\end{align*}

 Applying Lemma \ref{lemmamarkov} entails 
\begin{multline*}
\!\!\!\!\!\!\!  \pp\( \sup_{\gamma\in\cD_{m}\atop \cap[1+2\e,2-2\e]} \!\! \sup_{0\le v<u \le 1 \atop |v-u|\le 2^{-m}}  \!\! 2^{\frac{m}{ \gamma  (\be(\cM_{u+2^{-m}})+\frac 2m)}} \lba\int_v^u   \int_0^{2^{-\frac{m}{\gamma}}} z^{{1}/{\be(\cM_{s-})}} \wt N(ds,dz)\rba\ge 6m^2\) \\
\le C  \cdot 2^{m} \cdot e^{-m}.
\end{multline*}
Borel-Cantelli Lemma  yields that when $m$ becomes larger than some  $m_0$, for every $\gamma\in\cD_{m}  \cap[1+2\e,2-2\e] $ and $  |v-u|\le 2^{-m}$ (with $u>v$), 
$$\lba\int_v^u\int_0^{2^{-\frac{m}{\gamma}}} z^{{1}/{\be(\cM_{u-})}} \wt N(du,dz)\rba\le  6m^22^{\frac{-m}{ \gamma  (\be(\cM_{u+2^{-m}})+2/m)}} .$$
Assume that $r_0$ is so small that  $2r_0 \leq 2^{-m_0}$. By our choices for $n$ and $m$,   one has $m > n$, so  $\Upsilon^{n+1}_{J_{n+1}(t)} \in \mathcal{D}_m$. By choosing $\gamma = \Upsilon^{n}_{J_{n }(t)} $, $u=t+r$ and $v=t-r$, one gets
\begin{align*}
A_1(r) & \leq  6m^2  2^{\frac{-m}{ \Upsilon^{n}_{J_{n }(t)}  (\be(\cM_{t+r+2^{-m}})+2/m)}} \\
& \le  12  (\log_2(1/4r))^2  (2r)^{1/ ( \Upsilon^{n}_{J_{n }(t)}  (\be(\cM_{t+r+2^{-m}})+2/m))}.
\end{align*} 
In addition,  by continuity of $\cM$ at $t$, when $r_0$ is small enough, one has
 $$ |\Upsilon^{n}_{J_{n}(t)}  - \Upsilon(t)| < \ep_1 \   \ \mbox { and } \ \ \be(\cM_{t+r+2^{-m}})+2/m  \leq \be(\cM_{t}) +\ep_1,$$
  so finally 
 \begin{align*}
A_1(r) &  \leq \   r ^{\frac{1}{ (\Upsilon(t) + \e_1) (\be(\cM_{t}) +\ep_1)} }.
\end{align*}

On the other hand, recalling the constant $\e_0>0$ in Definition \ref{defstablelike}, an immediate computation shows that 
\begin{align*}
A_2(r) &\le \int_{t-r}^{t+r}  \int_0^{2^{-m/(\Upsilon(t) +\e_1) }}    z^{1/(\be(\cM_t)+\e_1)} \frac{dz}{z^2}\, ds  \\
&\le  \frac{2r}{1/\e_0 -1}  2^{\frac{-m}{(\Upsilon(t) +\e_1)}(\frac{1}{(\be(\cM_t)+\e_1)}- 1)}  \\
&\le \frac{1}{2/\e_0 -2}  (4r)^{\frac{1}{ (\Upsilon(t) + \e_1) (\be(\cM_{t}) +\ep_1)}  -  \frac{1}{\Upsilon(t) + \e_1} +1}\\
&\le \  r ^{\frac{1}{ (\Upsilon(t) + 2\e_1) (\be(\cM_{t}) +\ep_1)} },
\end{align*}
 as soon as $r_0$ is small enough. 

Combining these estimates, one obtains that for all $r\le r_0$, 
\begin{align*}
\cM_{t+r}-\cM_{t-r} & \le  r ^{\frac{1}{ (\Upsilon(t)+3\e_1 ) (\be(\cM_{t}) +\ep_1)}},
\end{align*}
which entails for all $0<\mathfrak{r}< \cM_{t+r_0}- \cM_{t-r_0}$,
\begin{align*}
\mu(B(\cM_t,\mathfrak{r})) \ge {\mathfrak{r}}^{(\Upsilon(t)+ 3\e_1)(\be(\cM_t)+\e_1)}.
\end{align*} 
One concludes that 
\begin{align*}
\ov \dim(\mu,\cM_t)&\le (\Upsilon(t)+3\e_1 )(\be(\cM_t)+\e_1 ).
\end{align*}
Letting $\e_1 \to 0 $ yields the desired upper bound for $\ov\dim(\mu,\cM_t)$.
\end{preuve}

\subsection{Dimension of $\cC(\Upsilon, \e')$ }\label{sectionbeforeextend}
Here we  prove that $\dimh \cC(\Upsilon, \e') \ge 2/(\Upsilon_{\min}+\e')-1$.

\begin{lem}\label{lemscaling} With probability $1$, for every $\Upsilon$ in Theorem \ref{theogeneral} with $\Upsilon_{\min}\in[1+2\e,2-2\e]$ and $\e'>0$, there exists a finite positive constant $K_{\Upsilon,\e'}$ such that for all $B\in\cB([0,1])$,  
\begin{align}\label{eqscaling}
\nu_{\Upsilon,\e'}(B)\le K_{\Upsilon,\e'} |B|^{\frac{2}{\Upsilon_{\min}+2\e'}-1}    
\end{align}
 \end{lem}

\begin{preuve} \ 
Let $\Upsilon$ and $\e'>0$ be fixed. 

Let $B$ be an open interval in $[0,1]$ such that  $|B| \le \eta_{n_0}$. 

If $B\cap \cC (\Upsilon,\e') = \emptyset$, \eqref{eqscaling} holds trivially. 

If $B\cap \cC (\Upsilon,\e') \neq \emptyset$, let $n_1$ be the largest integer such that $B$ intersects $ C_{n_1}(\Upsilon,\e') $ in exactly one basic interval,  denoted by ${J}^{n_1}$.

 Denote by $\de_{n_1+1}(\Upsilon,\e',{J}^{n_1})$ the minimal distance between any two intervals of $\cC_{n_1+1}(\Upsilon,\e')$ which are contained in  ${J}^{n_1}$. Then $|B|$ contains at most 
$$ \min\big(M_{n_1}(\Upsilon^{n_1}_{{J}^{n_1}}), |B|/ \de_{n_1+1}(\Upsilon,\e',{J}^{n_1})  \big)$$
 intervals of generation $n_1+1$.

 In addition, by construction, one has
 \begin{equation}
 \label{maj2}
 \de_{n_1+1}(\Upsilon,\e',{J}^{n_1})  \ge \frac{\eta_{n_1 }}{2M_{{n_1}}(\Upsilon^{n_1}_{{J}^{n_1}} )}.
 \end{equation}
 Hence by \eqref{defnu}, since all the intervals $J^{n_1+1}$ of generation $n_1+1$ within $J^{n_1}$ have the same $\nu$-mass, one has  (using \eqref{maj2})
\begin{align*}
&\nu_{\Upsilon,\e'}(B) \\
&\le \min\big(M_{n_1}(\Upsilon^{n_1}_{{J}^{n_1}}), |B|/ \de_{n_1+1}(\Upsilon,\e',{J}^{n_1})  \big)  \cdot  \nu_{n_1+1}(J^{n_1+1})\\
&\le \min\left( M_{n_1}(\Upsilon^{n_1}_{{J}^{n_1}}), |B|\frac{2M_{{n_1}}(\Upsilon^{n_1}_{{J}^{n_1}})}{\eta_{n_1,0}}   \right) \cdot  \( \prod_{k=n_0}^{n_1-1} \frac{1}{M_k(\Upsilon_{\min}+\e'  )} \) \cdot \frac{1}{M_{n_1}(\Upsilon^{n_1}_{{J}^{n_1}} )}   \\ 
&\le 2 \( \prod_{k=n_0}^{n_1-1} \frac{1}{M_k(\Upsilon_{\min}+\e'  )}\) \cdot  \eta_{n_1 }^{-1} \cdot \min \( \eta_{n_1 },   |B| \) .
\end{align*}
Due to our choices for the sequence $(\eta_{n})_{n\geq 1}$, when $n_0$ is large, 
$$  \prod_{k=n_0}^{n_1-1} \frac{1}{M_k(\Upsilon_{\min}+\e'  )}  \leq ({M_{n_1-1}(\Upsilon_{\min}+\e'  )}  )^{-1} \le \eta_{n_1 }^{\frac{2}{\Upsilon_{\min}+2\e'}-1}, $$
so applying the inequality  $x\wedge y \le x^s  y^{1-s}$ for $s\in(0,1)$ and $0<x,y<1$ yields
\begin{align*}
\nu_{\Upsilon,\e'}(B)  
 &\le 2 \eta_{n_1 }^{\frac{2}{\Upsilon_{\min}+2\e'}-1} \cdot \eta_{n_1 }^{-1} \cdot \eta_{n_1 }^{2-\frac{2}{\Upsilon_{\min}+2\e'}}  \cdot |B|^{\frac{2}{\Upsilon_{\min}+2\e'}-1} = 2 |B|^{\frac{2}{\Upsilon_{\min}+2\e'}-1}.
\end{align*}
\end{preuve}

\mk

Finally the mass distribution principle applied to the measure $\nu_{\Upsilon,\e'}$, which is supported by the Cantor set $\cC(\Upsilon,\e')$, allows one to conclude that 
$$\dimh \cC(\Upsilon,\e') \geq \frac{2}{\Upsilon_{\min}+2\e'}-1.$$

\subsection{Extension to $\Upsilon_{\min}\in\{1,2\}$}\label{sectiongamma=1,2}

 Letting $\e\to 0$ along a countable sequence yields that almost surely, for all $\Upsilon$ with $\Upsilon_{\min}\in(1,2)$,
\begin{align}\label{eqAtilde} \dimh F(\Upsilon) \ge \frac{2}{\Upsilon_{\min}} -1.
\end{align}


It remains to treat the extreme cases. 

\sk
{\bf First case : $\Upsilon_{\min}=1$.} For each $\e_2>0$, there exists an open interval $\cO\in(0,1)$ such that every $t\in\cO$ satisfies $ \Upsilon(t) \ge 1+ \e_2>1$. Applying \eqref{eqAtilde} yields that $\dimh F(\Upsilon) \ge \frac{2}{1+ \e_2}-1$ . Letting $\e_2\to 0$ establishes that  $\dimh F(\Upsilon) = 1$.

\sk
{\bf Second case : $\Upsilon_{\min}=2$,} i.e. $\Upsilon\equiv 2$. In order to prove $\dimh F(\Upsilon) \ge 0$, it suffices to show that there exists almost surely $t\in(0,1)$ such that $\ov \dim(\mu,\cM_t) = 2\be(\cM_t)$, i.e. $F(\Upsilon) \neq \emptyset$. To this end,   some changes are needed for the construction of the Cantor set.  We only sketch the proof since it is essentially the same as the one in the precedent sections (with simplification).  Set 
\begin{align*}
\begin{cases}
\rho_{0} =  1/2  \mbox{ and }  \rho_{n} = \exp(-\rho_{n-1}^{-1}) \mbox{ for all } n\ge 1,  \\
 \eta_{n} = \rho_n \log(1/\rho_n)^{-1} \mbox{ for all } n\ge 0. 
\end{cases}
\end{align*}
 Let $\cJ_n(2)$ be the set composed of intervals $J_{n,k} = [k\eta_n, (k+1)\eta_n)$  that  satisfy
\begin{align}\label{config11}
\begin{cases}
N\( [ {J}_{n, k-2}\times\[ \rho_n^{1/2}\( \log\frac{1}{\rho_n} \) ^{-3},1\] \) = 1, \\
N\( {J}_{n, k+2}\times\[ \rho_n^{1/2}\( \log\frac{1}{\rho_n} \) ^{-3},1\] \) = 1, \\
N\( \wh{J}_{n,k}\times\[ \rho_n^{1/2}\( \log\frac{1}{\rho_n} \) ^{-3},1\] \) =0, \\
\end{cases}
\end{align}
It is easy to check that any point  $t$ covered by the collection $\cJ_n(2)$ infinitely often satisfy $\ov \dim(\mu,\cM_t) \ge  2\be(\cM_t)$ (necessarily, one has equality thanks to Theorem \ref{theoexponent}). We construct as before,  by induction,  the collection $\cC_n(\Upsilon\equiv 2)$ of basic intervals and the Cantor set $\cC(\Upsilon\equiv 2)= \bigcap_n \cC_n(\Upsilon\equiv 2)$ contained in $F(\Upsilon)$. The same arguments as in Lemma \ref{chosen} give  a constant $C_n$ uniformly bounded below and above by $0$ and $+\infty$ such that for any fixed ${J_{n,k}}$,
\begin{align*}
\pp({J_{n,k}}\in\cJ_{n}(2)) =  C_n \cdot \rho_n \( \log \frac{1}{\rho_n}\) ^{4}.
\end{align*}
Thus one bounds from above the probability that there exists  ${J_{n,k}}$ such that none of the intervals $J_{n+1, k'}$ contained in ${J_{n,k}}$  belongs to $\cJ_{n+1}(2)$  by
\begin{align*}
\frac{1}{\eta_n} \(  1- C_{n+1}\cdot\rho_{n+1} \( \log \frac{1}{\rho_{n+1}}\) ^{4}  \) ^{\frac{\eta_n}{\eta_{n+1}}}.
\end{align*}
Observe that 
$$\frac{\eta_n}{\eta_{n+1}} = C_{n+1} \rho_{n-1}^{-3}\cdot (C_{n+1})^{-1}\rho_{n+1}^{-1} \( \log\frac{1}{\rho_{n+1}} \) ^{-3} \mbox{ with } C_{n+1}\rho_{n-1}^{-3} \gg 1.  $$
 So the probability in question is less than 
$$\eta_n^{-1} e^{- C_{n+1} \rho_{n-1}^{-3}} \le \eta_n^{-1} e^{- 3\rho_{n-1}^{-1}} =  \eta_n^{-1} \rho_n^3 \le \rho_n.$$
Borel-Cantelli Lemma implies the existence of a sequence of embedded interval with length tending to $0$ that satisfy \eqref{config11}. This justifies that $F(\Upsilon)\neq\emptyset$.


\section{Space spectrum : proof of Theorem \ref{theospace} } 

\label{sec_space}

\subsection{A first theorem on dimensions, and the space spectrum}

  Throughout this section, we set $\e= \e_0$, which is defined in \eqref{beta}. 
 We are going to prove the following theorem.

\begin{theo}
\label{th_dimension}
Let $\ep>0$.  Denote by $\cP = \{(T_n,Z_n)\}_{n\geq 1}$  a Poisson point process that generates the Poisson measure $N(dt,dz)$  with intensity $dt\otimes dz/z^2$.  Consider the family \eqref{defLalpha} of stable processes $(\cL^\al_.)_{\alpha\in (\ep,1-\ep)}$. Also,  for every  non decreasing 
c\`adl\`ag  function $f:[0,1] \to [\ep,1-\ep]$, consider the process 
\begin{equation}
\label{defLf}
\cL^f_t = \int_0^t\int_0^1 z^{\frac{1}{f(t-)}} N(ds,dz).
\end{equation}

Almost surely, for every set $E\subset [0,1]$, for every function $f:[0,1] \to[\ep,1-\ep]$, if $\alpha  <  \inf\{f(t): t\in E\}$ and $\beta >   \sup\{f(t): t\in E\}$, then
$$ \dimh(\cL^\alpha(E) )\leq \dimh \cL^f(E) \leq \dimh(\cL^\beta(E)) = \frac{\beta}{\alpha}   \dimh(\cL^\alpha(E) ).$$

\end{theo}

Before proving Theorem \ref{th_dimension}  next subsection, let us explain how we deduce the space spectrum of the occupation measure.

\mk

As mentioned in the introduction, almost surely, for every set $E\subset\R$, the image of $E$ by an $\alpha$-stable process $\cL^\alpha$ has Hausdorff dimension $\alpha \dimh (E)$.
Applying Theorem \ref{th_dimension} to the function $f(t) = \beta(\cM_t)$, which is almost surely c\`adl\`ag, one is now ready to prove Theorem \ref{th_dim_images}.

\mk

\begin{preuve}{\it \ of Theorem \ref{th_dim_images} : } \ 
The first part (the formula \eqref{eqfinal}) is immediate.

For the second part, let $E\subset \zu$ be such that for every non-trivial subinterval $I\subset [0,1]$, $\dimh (E) = \dimh (E\cap I)$.

For every $\eta>0$, there exists an interval $I$ of length less than $\eta$ such that 
$$\Big[\inf_{t\in E\cap I } \beta(\cM(t)),\sup_{t\in E\cap I} \beta(\cM(t-)) \Big] \subset \Big [\sup_{t\in E} \beta(\cM(t-)) -\eta,\sup_{t\in E} \beta(\cM(t-)) \Big] .$$
This follows from the c\`adl\`ag regularity of $t\mapsto\beta(\cM(t))$. Hence, applying \eqref{eqfinal} to $E\cap I$ gives
 $$\dimh \cM(E\cap I)   \in \dimh(E\cap I)  \cdot  \Big [\sup_{t\in E} \beta(\cM(t-)) -\eta,\sup_{t\in E} \beta(\cM(t-)) \Big]  .$$
Since $\dimh \cM(E\cap I) \leq \dimh \cM(E) $ and  $\dimh (E) = \dimh (E\cap I)$, the result follows by letting $\eta$ tend to zero.
\end{preuve}

\sk

One deduces a corollary from Theorem \ref{th_dim_images}, which will be used in the proof of Theorem \ref{theospace}.

\begin{cor}\label{cor2}  For every open interval $I \subset [0,1]$ and  $h\ge 0$,  consider the smallest interval $I_0\subset I$ (it may be not open or reduced to a point) such that  $\ov E^t_\mu(I, h) = \ov E^t_\mu(I_0, h)$.  Denote by $d(I_0)$ the right endpoint of $I_0$
 Almost surely,  for every open interval $I \subset [0,1]$ and  $h\ge 0$, 
one has
 \begin{equation}\label{eqcor2}
 \dimh \cM\( \ov E^t_\mu(I, h) \) =  \dimh \ov E^t_\mu(I_0, h)   \cdot  \sup_{t\in I_0 \setminus {d(I_0)}} \be(\cM_t).
 \end{equation}
 with the convention that $0\times (-\infty) =0 $ and $ (-\infty)\times (-\infty) = -\infty$.
 \end{cor}

\begin{preuve} \   If $\ov E^t_\mu(I, h)$ is empty or a singleton, there is nothing to prove.  One thus assumes that 
$\ov E^t_\mu(I, h)$ is neither empty nor a singleton, so  $I_0$ is a non-trivial interval.  One could check the analysis in the proof of Theorem \ref{theotime} for a construction of $I_0$.   Observe that the left-hand side of \eqref{eqcor2} is smaller than the right-hand side due to Theorem \ref{th_dim_images}.  The converse inequality follows by minimality of  $I_0$  and a localization procedure as in the proof of Theorem \ref{th_dim_images}.
\end{preuve}

\mk

\begin{preuve}{\it \ of Theorem \ref{theospace} and \ref{theospaceextreme} : }  To deduce the space spectrum, one needs some additional analysis other than the time spectrum.  This is due to the following basic observation  :  for all $t\in S(\cM)$, $\cM_{t-}$ is not in the range of $\cM$, but in the support of $\mu$.  

When $\cO$ does not intersect the range of $\cM$,  the level set $\ov E_\mu(\cO,h)=\emptyset$, as is given in Theorem \ref{theospace} and \ref{theospaceextreme}. 

When $\cO$ intersects the range of $\cM$, by the c\`adl\`ag property of $\cM$,  there is a non-trivial interval $\wt\cO$ such that 
$$\cM((0,1))\cap \cO = \cM(\wt\cO).$$
The set $\wt\cO$ is  an open set  $(a,b)$ if $\cM$ enters $\cO$ continuously, or is a semi-open interval $[a,b)$ if $\cM$ enters  $\cO$ with a jump.  In any case, $\cM_{a-} \notin \cO$ and $\cM_b\notin \cO$ because $\cO$ is open.  
Observe that 
\begin{align*}
&\ov E_\mu(\cO, h)  \\
&= \{x\in \supp\mu \cap \cO : \ov \dim(\mu, x) = h \} \\
&= \{\cM_t\in \cO : \ov\dim(\mu,\cM_t) = h \} \cup \{\cM_{t-}\in\cO : t\in S(\cM) \mbox{ and } \ov\dim(\mu,\cM_{t-}) = h  \} \\
&= \cM(\ov E^t_\mu(\wt\cO,h)) \cup \{ \cM_t : t\in S(\cM)\cap \wt\cO \mbox{ and } \ov\dim(\mu, \cM_{t-}) = h\} 
\end{align*}
Since the Hausdorff dimension of the second set in the last union is at most $0$,  one distinguishes two types of situations according to the value of $h$. 

\noindent $\bullet$ Type A.  The time level set $\ov E^t_\mu(\wt\cO, h)\neq\emptyset$, so one ignores the second set in the last union when computing the Hausdorff dimension of $\ov E_\mu(\cO, h)$.   If $\# \ov E^t_\mu(\wt\cO, h)=1$ (necessarily $h=\be(\cM_t)$ for some $t\in S(\cM)$ with $\be(\cM_t)\ge 2\be(\cM_{t-})$), one has $\ov d_\mu(\cO, h)=0$ which  coincide with the formula in Theorem \ref{theospaceextreme}.    Otherwise Corollary \ref{cor2} applied to $h$ and $\wt \cO$ entails the existence of a minimal $\wt\cO_0$ (that we can and will suppose open)  such that 
\begin{align*}
\ov d_\mu(\cO, h) &= \ov d^t_\mu(\wt\cO_0,h) \cdot \sup_{t\in\wt\cO_0\setminus d(\wt\cO_0)} \be(\cM_t) \\
& =  \sup\llb  \hat{g}_\al(h)/\al :  \al\in \{\be(\cM_t) : t\in \wt \cO_0 \}  \rrb \cdot \sup_{t\in\wt\cO_0} \be(\cM_t) \\
&= \sup\llb  \hat{g}_\al(h) : \al \in \{ \be(\cM_t) : t\in \wt\cO_0 \}\rrb \\
& =  \sup\llb  \hat{g}_\al(h) : \al \in \{ \be(\cM_t) : t\in \wt\cO \}\rrb \\
& =  \sup\llb  \hat{g}_\al(h) : \al \in \{ \be(\cM_t) : \cM_t\in \cO \}\rrb,
\end{align*}
as desired.

\noindent $\bullet$ Type B. The time level set 
\begin{equation} \label{eqtypeb}\ov E^t_\mu(\wt\cO, h) =\emptyset  \end{equation} so one has to  consider the (at most) countable set  $\cH = \{ \cM_t : t\in S(\cM)\cap \wt\cO \mbox{ and } \ov\dim(\mu, \cM_{t-}) = h\} $.   Compared with the time spectrum,   several cases may occur according to the value of $h$.  Recall that  $\wt\cO= [a,b)$ (when $\cM$ jumps into $\cO$) or $(a,b)$  (when $\cM$ enters $\cO$ continuously).  In the following analysis,  the case $\cM_{t_0} \notin\cO$ is trivial.  We thus assume that every $\cM_{t_0}$ below belongs to $\cO$. 

\begin{enumerate}
\item \eqref{eqtypeb} is due to $2\be(\cM_{t_0-})<h<\be(\cM_{t_0})$ with $t_0\in S(\cM)$.   For all $t>t_0$,  $\ov \dim(\mu, \cM_{t-})\ge \be(\cM_{t-})> \be(\cM_{t_0})>h$.   For all $t\le t_0$,  $\ov \dim(\mu, M_{t-}) \le 2\be(\cM_{t-}) \le 2\be(\cM_{t_0-})<h$.   So $\cH = \emptyset$ and $\ov d_\mu(\cO, h) = -\infty$ as desired.  

\sk

\item \eqref{eqtypeb} is due to the fact that $2\be(\cM_{t_0-})< \be(\cM_{t_0})=h$ with $t_0\in S(\cM)$ and $\ov\dim(\mu, \cM_{t_0}) \neq \be(\cM_{t_0})$.  As in the last item,  $\cH = \emptyset$ as desired. 

\sk 

\item \eqref{eqtypeb} is due to  $h= 2\be(\cM_{t_0-}) < \be(\cM_{t_0})$ or $h= 2\be(\cM_{t_0-}) = \be(\cM_{t_0})$.  As before for all $t\neq t_0$,  one has $\ov\dim(\mu, \cM_{t-}) \neq h$.   If  $\ov\dim(\mu, \cM_{t_0-}) = 2\be(\cM_{t_0-})$,  $\cH = \{t_0\}$,  otherwise $\cH = \emptyset$.  This coincides with  Theorem \ref{theospaceextreme}.  

\sk

\item \eqref{eqtypeb} is due to  $h\ge 2\be(\cM_{b-})$. For all $t<b$,  $\ov \dim(\mu, \cM_{t-})\le 2\be(\cM_{t_0-}) < 2\be(\cM_{b-})\le h$, hence $\cH= \emptyset$, as desired.

\sk

\item  \eqref{eqtypeb} is due to  $h\le \be(\cM_a)$.  Recall first that $\be(\cM_a)\notin \cO$. Further, for all $t>a$, $\ov\dim(\mu, \cM_t) > h$. Hence,  $\cH=\emptyset$, as desired. 
\end{enumerate}
  
\end{preuve}


\subsection{Proof of Theorem \ref{th_dimension}}

We start with a Lemma describing the distribution properties of the Poisson point process $\mathcal{P}=\{(T_n,Z_n)\}_{n\geq 1}$.

\begin{lem}
\label{lem_distrib}
For every $j\geq 1$, let $\mathcal{P}_j = \{n: Z_n\in [2^{-j-1},2^{-j})\}$. Almost surely, there exist two  positive decreasing sequences $(\ep_j)_{j\geq 1}$ and $(\eta_j)_{j\geq 1}$  converging to zero such that for every integer $J$ large enough, one has:
\begin{enumerate}
\item
  $ 2^{J(1-\ep_J)} \leq \#  \mathcal{P}_J  \leq 2^{J(1+\ep_J)} $,

\item
for  every interval $I\subset [0,1]$  with length $2^{-J
}$, 
$$ 1  \leq \#  \bigcup_{ j\leq J(1+\eta_J)}   \{n\in \mathcal{P}_j : T_n\in I\}  \leq 2^{J\ep_J} ,$$

\item
for   every interval $I\subset [0,1]$  with length $2^{-J
}$, 
$$ 0\leq   \#  \bigcup_{ j\leq J/3}   \{n\in \mathcal{P}_j : T_n\in I\}  \leq 1 ,$$

\item
for every interval $I\subset [0,1]$  with length $2^{-J }$, for every $ j\geq J(1+\eta_J)$,
$$   \#     \{n\in \mathcal{P}_j : T_n\in I\}  \leq 2^{j(1+\ep_j)}2^{-J  
}.$$

\end{enumerate}

\end{lem}
 
Routine computations as in Lemma \ref{chosen} entail  Lemma \ref{lem_distrib}.
 

 \mk

Let $E\subset \zu$, $\alpha = \inf\{f(t): t\in E\}$ and $\beta=  \sup\{f(t): t\in E\}$.

\mk

Call $E^\alpha $ (resp. $E^f$, $E^\beta$) the image of $E$ by $\cL^\alpha$ (resp. 
$\cL^f$, $\cL^\beta$).

\mk

\begin{lem}
\label{lem1}
Almost surely, the following holds. 
With each interval $\widetilde B^\alpha $ such that $\widetilde B^\alpha\cap E^\alpha \neq \emptyset$, one can associate  an interval of the form $B^\alpha=\cL^\alpha([T_m, T_n))$ such that $|B^\alpha| \leq 2|\widetilde B^\alpha |$ and possibly a singleton of the form $\{ \cL^\alpha(T_n)\}$, such that 
$$ E^\alpha \cap (B^\alpha \cup \{ \cL^\alpha(T_n)\}) = E^\alpha\cap \widetilde B^\alpha.$$
The same holds true for every interval $\widetilde B^f$ such that $E^j\cap \widetilde B^f \neq\emptyset$, which can be replaced by $B^f=\cL^f([T_m, T_n))$ and   possibly a singleton.
\end{lem}

\begin{preuve}
 
Almost surely all the processes $\cL^\alpha$, $\cL^\beta$ and $\cL^f$ are strictly increasing and c\`adl\`ag.

Let $\widetilde  B^\alpha =[x^\alpha,y^\alpha] $ be   an interval satisfying  $\widetilde B^\alpha\cap E^\alpha \neq \emptyset$.

If $x^\alpha$ is not of the form $\cL^\alpha(T_m)$, then two cases occur:
\begin{itemize}
\item
when $x^\alpha \notin \cL^\alpha(E)$:   $\widetilde  B^\alpha$ can be replaced by $[x'^\alpha,y^\alpha]$, where $x'^\alpha = \inf (\widetilde  B^\alpha\cap E^\alpha)$, without altering the covering $\cR^\alpha$. Since $\cL^\alpha$ is increasing and c\`adl\`ag, $x'^\alpha$ is necessarily the image of some jump point $T_m$ by $\cL^\alpha$.

\item
when $x^\alpha \in \cL^\alpha(E)$:   $x^\alpha$ can be written as $\cL^\alpha(t)$, for some $t$ which is a continuous time for $\cL^\alpha$.   Using the density of the jump points, there exists $(T_m,Z_m)$ such that $ T_m<t$  and $\cL^\alpha(t) - \cL^\alpha(T_m) <|\widetilde  B^\alpha|/2$. We then choose $x'^\alpha = \cL^\alpha(T_m)$.
\end{itemize}

 In all cases, $\widetilde   B^\alpha  $ is replaced by $B'^\alpha=[x'^\alpha,y^\alpha]$, where $| B'^\alpha| \leq3/2  |B^\alpha| $.

Similarly, if $y^\alpha$ is not of the form $\cL^\alpha(T_n-)$ (i.e. the left limit of $\cL^\alpha$ at $T_n$ for some jump point $T_n$),  then :
\begin{itemize}
\item
when $y^\alpha \notin \cL^\alpha(E)$:   $B'^\alpha$ can be replaced by $B^\alpha=[x'^\alpha,y'^\alpha]$, where $y'^\alpha = \sup (B^\alpha\cap E^\alpha)$, without altering the covering $\cR^\alpha$. Since $\cL^\alpha$ is increasing and c\`adl\`ag, $y'^\alpha$ is of the form  $\cL^\alpha(T_n-)$ for some jump point $T_n$.

\item
when $y^\alpha = \cL^\alpha(T_n) $ for some jump time $T_n$:  Then $B'^\alpha$ can be replaced by $\{\cL^\alpha(T_n)\} \cup B^\alpha$, where $B^\alpha=[x'^\alpha,\cL^\alpha(T_n-)] $. Indeed, there is no point of $E^\alpha $ between $\cL^\alpha(T_n-)$ and $\cL^\alpha(T_n)$.

\item
when $y=\cL^\alpha(t)$ for some $t$ which is a continuous time for $\cL^\alpha$.   Using the same argument as above,  there exists $(T_n,Z_n)$ such that $ T_n>t$  and $\cL^\alpha(T_n) - \cL^\alpha(t) <|B^\alpha|/2$. We then choose $y'^\alpha = \cL^\alpha(T_n-)$.
\end{itemize}

This proves the claim. 
\end{preuve}

Observe that the previous Lemma holds almost surely, for every interval $B^\alpha$, for all $\alpha$, since the randomness is only located in the distribution of the Point Poisson process and the  strictly increasing and c\`adl\`ag properties of the processes, which hold simultaneously   almost surely.

\mk
Next Lemma establishes that the increment of the process in an interval $I$ is approximately the same order as the size of the largest jump in $I$, uniformly for all $I$ and all the parameters. 

\begin{lem}
\label{lempoisson}
With probability one, there exists a non-decreasing function $g: \zu\to \R^+$ with $g(0)=0$, continuous at 0,   such that the following holds. Let $(T_m,Z_m)$ and $(T_n,Z_n)$ (with $T_m<T_n$) be two couples of the point Poisson process. Let  $B^\alpha=[\cL^\alpha(T_m),\cL^\alpha(T_n-)]$, $B^\beta=[\cL^\beta(T_m),\cL^\beta(T_n-)]$ and $B^f=[\cL^f(T_m),\cL^f(T_n-)]$. Then when $B^\alpha=[\cL^\alpha(T_m),\cL^\alpha(T_n-)]$ is small enough,  one has
\begin{equation}
\label{inegbalpha}
  |B^\alpha|^{\alpha/\beta+ g ( |B^\alpha|)}   \leq   |B^\beta|  \leq  |B^\alpha|^{\alpha/\beta-g( |B^\alpha|)}
  \end{equation}
and
\begin{equation}
\label{inegbalpha2}
 |B^\alpha| \leq  |B^f|  \leq   |B^\beta| .
\end{equation}

\end{lem}

\begin{preuve}

The three processes $\cL^\alpha$, $\cL^f$ and $\cL^\beta$ are almost surely pure jump processes with finite variations. One deduces that
\begin{equation}
\label{eq14}
|B^\alpha| = \sum_{p\in \N:T_p\in [T_m,T_n)} Z_p^{1/\alpha},
\end{equation}
and the same  holds true for $|B^\beta|$ by replacing $1/\alpha$ by $1/\beta$. Similarly,
$$|B^f| = \sum_{p\in \N:T_p\in [T_m,T_n)} Z_p^{1/f(T_p-)}.$$

Then \eqref{inegbalpha2} follows immediately since $f$ is monotone and  $\alpha \leq f(t) \leq \beta$.

\mk

We write $B=[T_m,T_n)$, and consider  $J$  the unique integer such that $2^{-(J+1)} < |T_n-T_m| \leq 2^{-J}$.  We assume that $J$ is so large that  $\ep_J\leq (1/(1-\ep)-1)/4 \ \leq (1/\alpha-1)/4$.

We now make use of Lemma \ref{lem_distrib}.

Let $(T_N,Z_N)$ be the point Poisson process in the above sum \eqref{eq14} with largest jump $Z_N$. We write $Z_N = 2^{-J_N}$. Then one decomposes $|B^\alpha|$ into
\begin{equation}
\label{decompBalpha}
|B^\alpha| = Z_N^{1/\alpha} + \sum_{ j\leq J(1+\eta_J):  T_p\in B \mbox{ and } p\in \mathcal{P}_j} Z_p^{1/\alpha} + \sum_{ j> J(1+\eta_J):  T_p\in B \mbox{ and } p\in \mathcal{P}_j} Z_p^{1/\alpha}.
\end{equation}

Call $S_1$ and $S_2$ the two above sums.

\mk

Assume that $J_N <J/3$.

Observe that  since $B$ strictly contains an interval of length $2^{-J-1}$, the left inequality part (2) of Lemma \ref{lem_distrib}  yields that $J_N \leq (J+1)(1+\eta_{J+1})$. 

Since $B$ is contained in an interval of length $2^{-J}$, one knows that  all the   jumps other than $(T_N,Z_N)$  appearing in formula \eqref{decompBalpha} are smaller than $2^{-J/(3\alpha)}$. Hence, the right inequality in part (2) of Lemma \ref{lem_distrib} yields
$$S_1 \leq 2^{J\ep_J}   2^{-J/(3\alpha)}. $$ 
Similarly,  applying part (4) of Lemma \ref{lem_distrib}, 
 $S_2$ is bounded by
\begin{eqnarray}
\nonumber S_2&  \leq   & \sum_{ j\geq J(1+\eta_J)} 2^{j(1+\ep_j)}2^{-J  } 2^{-j/\alpha}  \leq \frac{2^{J(1+\eta_J)(1+\ep_J-1/\alpha)-1}}{2^{3/4(1/\alpha-1)}-1}\\
\label{eq112}
  & \leq &  C_\ep 2^{-J /\alpha } 2^{  J(  \ep_J+ \eta_J(1+\ep_J-1/\alpha))},
  \end{eqnarray}
where $C_\ep:= \frac{1}{2^{3/4(1/(1-\ep)-1)}-1}$.      Recalling that $J_N\leq J/3$, one gets
\begin{eqnarray*}
    |B^\alpha|   & \leq  &  (Z_N)^{1/\alpha}  +2^{-J(1/(3\alpha)-\ep_J)} +C_\ep 2^{-J (1 /\alpha - \ep_J-  \eta_J(1+\ep_J-1/\alpha))} 
   \\
   & \leq &  (Z_N)^{1/\alpha}  +(Z_N)^{ 1/(\alpha)-\ep_{J_N}/3} +  (Z_N)^{  3(1 /\alpha - \ep_{J_N}- \eta_{J_N}(1+\ep_{J_N}-1/(1-\ep)))}. 
   \end{eqnarray*}
One concludes that 
 $$
(Z_N)^{1/\alpha} \leq    |B^\alpha|      \leq      (Z_N)^{1/\alpha -\widetilde\ep_{J_N} }  ,$$
 for some $\widetilde\ep_{J_N} $ which depends only on $\ep_{J_N} $ and $\eta_{J_N}$ (not on $\alpha$),  is decreasing as a function of $\ep_{J_N}$ and $\eta_{J_N}$, and which tends to zero  when $J_N$ tends to infinity. In addition, the fact that  $(Z_N)^{1/\alpha} \leq    |B^\alpha|   $ implies that $ {J_N} \geq  - \alpha \log_2 |B^\alpha|  \geq   \lfloor - \ep \log_2 |B^\alpha| \rfloor $. Hence  $\widetilde\ep_{J_N}  \leq  g_1(  \lfloor - \ep \log_2 |B^\alpha| \rfloor )$,  where $g_1 (r) = \widetilde  \ep_{ \lfloor- \ep \log_2  r \rfloor }$. One can write finally
\begin{equation}
\label{eq15}
(Z_N)^{1/\alpha} \leq    |B^\alpha|      \leq      (Z_N)^{1/\alpha - g_1 (|B^\alpha| )}  .
\end{equation}
 By construction, this mapping $g_1$ is non decreasing with $r$, and tends to 0 when $r$ tends to 0.

\mk

Observe that since $\widetilde \ep_{J_N}$ is small (uniformly in $\alpha$), one also has 
\begin{equation}
\label{eq16}
(Z_N)^{1/\alpha} \leq    |B^\alpha|      \leq      (Z_N)^{1/(2\alpha )}  .
\end{equation}

\mk
Assume now  that $J_N \geq J/3$.

All the   jumps other than $(T_N,Z_N)$  involved in formula \eqref{decompBalpha} are smaller than $2^{-J_N/\alpha}$. Hence, part (2) of Lemma \ref{lem_distrib} yields
$$S_1 \leq 2^{J\ep_J}   2^{-J_N/\alpha} \le 2^{-J_N(1/\al-3\e_{J_N})}. $$ 
The sum
 $S_2$ is still bounded by above by \eqref{eq112} with $J$ replaced by $3J_N$. One deduces that 
\begin{eqnarray*}
    |B^\alpha|  
    \leq &  (Z_N)^{1/\alpha}  +(Z_N)^{ 1/(\alpha) - 3\ep_{J_N}} +  (Z_N)^{  3(1 /\alpha - \ep_{J_N}- \eta_{J_N}(1+\ep_{J_N}-1/(1-\ep)))}. 
   \end{eqnarray*}
One concludes that 
 $
   |B^\alpha|      \leq      (Z_N)^{1/\alpha -\widetilde\ep _{J_N}}  $ for some $\widetilde\ep _{J_N} $ which depends only on $J_N$ (not on $\alpha$), and which tends to zero when $J_N$ tends to infinity. For the same reasons as above, equation \eqref{eq15} holds true.

\mk

Using that \eqref{eq15} holds true with $\beta$ instead of $\alpha$ (but with the same mapping $\widetilde\ep $), one sees that 
$$
  |B^\beta| \leq     (Z_N)^{1/\beta - g_1 (|B^\beta| )}\leq   |B^\alpha|^{\alpha/\beta - \alpha g_1 (|B^\beta| )} .
$$
 In addition,  using \eqref{eq16} with $\beta$ instead of $\alpha$, one has $|B^\beta|  \leq (Z_N)^{1/2\beta} \leq |B^\alpha|^{\alpha/(2\beta)}$. We deduce that $\alpha g_1 (|B^\beta| )  \leq  (1-\ep)g_1(|B^\alpha|^{\alpha/(2\beta)})  := g_2( |B^\alpha|)$, hence
\begin{equation}
\label{eq11}
  |B^\beta| \leq     |B^\alpha|^{\alpha/\beta - (1- \ep)   g_2 (|B^\alpha| )}  .
\end{equation}

Similarly, recalling that $|B^\alpha|$ and $|B^\beta|$ are small quantities, 
\begin{equation}
\label{eq11'}
  |B^\beta| \geq     (Z_N)^{1/\beta } \geq   |B^\alpha|^{1/( \beta (1/\alpha-   g_1 (|B^\alpha| ))}\geq   |B^\alpha|^{\alpha/\beta + 2\beta g_1(|B^\alpha| )} \geq   |B^\alpha|^{\alpha/\beta +  g_3(|B^\alpha| )}  .
\end{equation}
where $g_3(r) = 2(1-\ep)  \widetilde\ep (r )$.
Finally, \eqref{eq11} and \eqref{eq11'} gives the result, with $g(r) =  \max(g_2(r),g_3(r))$.
\end{preuve}

Observe that one can also write 
\begin{equation}
\label{inegbalpha3}
  |B^\alpha|^{\alpha/\beta+ \widetilde g ( |B^\beta|)}   \leq   |B^\beta|  \leq  |B^\alpha|^{\alpha/\beta-\widetilde g( |B^\beta|)}
  \end{equation}
for some mapping $\widetilde g$ which enjoys the same properties as $g$.

\mk

One can now prove Theorem \ref{th_dimension}.

\mk

The following holds almost surely, since it depends only on Lemmas \ref{lem1} and \ref{lempoisson}.
\mk

Let us denote by $d_\alpha = \dimh E^\alpha$, $d_\beta = \dimh E^\beta$, $d_f = \dimh E^f$.

\sk

Let $s>  d_\alpha \,\beta/\alpha    $, and let $\widetilde s=  s  \alpha/\beta -   (s \alpha   /\beta-d_\alpha )/2$. One has $  d_\alpha   < \widetilde s < s \alpha/\beta$.

\sk

By definition of $d_\alpha$, there exists $\eta>0$ such that $\mathcal{H}^{\widetilde s}_{\eta/2} (E^\alpha) \leq 4^{-s} .$  Hence, for some $\eta/2$-covering 
$  \cR^\alpha$ of $E^\alpha$, one has
$$\sum_{  \widetilde B^\alpha \in   \cR^\alpha} |   \widetilde B^\alpha|^{\widetilde s} \leq 2^{-s}.$$
 
First,  using \ref{lem1}, by slightly modifying the intervals $\widetilde B^\alpha \in \widetilde \cR^\alpha$,   one can replace these intervals with intervals of the form $B^\alpha=\cL^\alpha([T_m, T_n))$ (plus at most a countable number of singletons), satisfying $|B^\alpha| \leq 2 |\widetilde B^\alpha|$, whose union is still covering $E^\alpha$. 

Hence, the initial $\eta/2$-covering $\widetilde \cR^\alpha$ can be replaced by an $\eta$-covering  $\cR^\alpha$, such that  
  one has 
$$\sum_{   B^\alpha \in   \cR^\alpha} |  B^\alpha|^{\widetilde s} \leq 1.$$

Let us choose $\eta$ so small that $g  (\eta^{\alpha/\beta-g(\eta)}) < \frac{s \alpha   /\beta-d_\alpha }{2 s}$.

Each ball $B^\alpha$ is written $\cL^\alpha(B)$, where $B=[T_m,T_n)$. As above, we write   $B^\beta = \cL^\beta(B)$ and  $B^f = \cL^f(B)$, and  \eqref{inegbalpha}  and  \eqref{inegbalpha2} hold true.

Since the balls $(B^\alpha)$ form an $\eta$-covering of $E^\alpha$,  the balls $(B^\beta)$ form a $\widetilde \eta:=\eta^{\alpha/\beta-g(\eta)} $-covering of $E^\beta$, and the balls $(B^f)$ also form a $\widetilde \eta$-covering of $E^f$. We denote by $\mathcal{R}^\beta$ and $\mathcal{R}^f$ these two coverings.  One has
\begin{eqnarray*}
\sum_{  B^f \in   \cR^f} |  B^f|^{ s}  & \leq &  \sum_{  B^\beta \in   \cR^\beta} |  B^\beta|^{ s}   \leq \sum_{  B^\alpha \in   \cR^\alpha}  |  B^\alpha|^{ s (\alpha/\beta -g(|B^\alpha|))}\\
& \leq&  \sum_{  B^\alpha \in   \cR^\alpha}  |  B^\alpha|^{ s  \alpha/\beta -   (s \alpha   /\beta-d_\alpha )/2 } = \sum_{  B^\alpha \in   \cR^\alpha}  |  B^\alpha|^{ \widetilde s} \leq 1 .   
\end{eqnarray*}
 
Since $\cR^\beta$ is  an $\widetilde \eta $-covering of $E^\beta$, the  $s$-pre-Hausdorff measure of  $E^\beta$, $\mathcal{H}^s_{\widetilde \eta}(E^\beta)$ is less than $1$.  The same holds for $\mathcal{H}^s_{\widetilde \eta}(E^f)$.
This remains true for any sufficiently small $\widetilde \eta>0$, we conclude that both $\mathcal{H}^{s} (E^f)$ and $ \mathcal{H}^{s} (E^\beta) $ less than 1, hence  $d_f$ and $d_\beta $ are smaller than $ s$. Since this holds   for any  $s>d_\alpha\,\beta/\alpha$, one gets that $ \max(d_f, d_\beta)\leq d_\alpha \,\beta/\alpha $.

\mk

Next, starting with a $\eta$-covering of $E^f$ by balls $B^f$, one associates with every ball $B^f = \cL^f([T_m,T_n))$  the ball $B^\beta = \cL^\beta([T_m,T_n))$, the same lines of computation (simply using that $|B^f|\leq |B^\beta|$) yields that $d_f\leq d_\beta$. 

The same argument shows that $d_\alpha\leq d_f$.

\mk

It remains us to prove the last inequality $ d_\alpha \leq d_\beta \,\alpha/\beta$.   The proof follows exactly the same lines, we write it without details.

Let $s>  d_\beta \,\alpha/\beta    $, and let $\widetilde s=  s \beta/\alpha -   (s\beta  /\alpha-d_\beta )/2$. One has $  d_\beta   < \widetilde s < s \ \beta/\alpha$.

\sk

There exists an  $\eta$-covering $\mathcal{R}^\beta$ of $E^\beta$ by intervals of the form $B^\beta = \cL^\beta([T_m,T_n))$, such that
$$\sum_{   B^\beta \in   \cR^\beta} |  B^\beta|^{\widetilde s} \leq 1.$$
One considers the associated intervals $(B^\alpha)$ and $(B^f)$, and the natural coverings $\mathcal{R}^\alpha$ and $\mathcal{R}^f$ of $E^\alpha$ and $E^f$ provided by these intervals.

 Let   $\eta$ be so small that $\widetilde g(|B^\beta|) \leq \widetilde g  (\eta ) < \frac{s -d_\beta \alpha/\beta }{4 s}$, where $\widetilde g$ is given by \eqref{inegbalpha3}.  One has
\begin{eqnarray*}
  \sum_{  B^\alpha \in   \cR^\alpha} |  B^\alpha|^{ s}     & \leq & \sum_{  B^\beta \in   \cR^\beta}  |  B^\beta|^{ s /(\alpha/\beta +\widetilde g(|B^\beta|))} \leq \sum_{  B^\beta \in   \cR^\beta}  |  B^\beta|^{ s  \beta/\alpha  -2 s \widetilde g(|B^\beta|)  \beta/\alpha } \\
& \leq&  \sum_{  B^\beta \in   \cR^\beta}  |  B^\beta|^{ s  \beta/\alpha  -  (s \beta/\alpha -d_\beta)/2 } = \sum_{  B^\beta \in   \cR^\beta}  |  B^\alpha|^{ \widetilde s} \leq 1 .   
\end{eqnarray*}

This holds true for any $\eta>0$ small enough, so that $\mathcal{H}^{s} (E^\alpha)   <+\infty$, hence  $d_\alpha \leq s$. Since this holds true for any  $s>  d_\beta \,\alpha/\beta   $, one gets that $d_\alpha \leq d_\beta \,\alpha/\beta  $.

 \mk
 
 One concludes that $d_\alpha\leq d_f \leq d_\beta = d_\alpha\beta/\alpha$.

\begin{remark}
It is certainly possible to short-cut the end of the proof, since one knows that $\alpha$ being fixed, almost surely, for every set $E\subset \zu$, $\dimh \cL^\alpha(E) = \alpha\dimh (E)$.

But since we consider all $\alpha$'s, it was easier to prove all inequalities at once.

\end{remark}



\bibliographystyle{spmpsci}   
\bibliography{xyangbiblio}



\end{document}